\documentclass[12pt,leqno]{article}
\usepackage{amssymb}
\usepackage[mathscr]{eucal}
\usepackage{amsmath,amssymb,latexsym,theorem,bbm}
\usepackage{color}
\usepackage{appendix}

\newcommand{\NN}{\mathbb{N}}
\newcommand{\RR}{\mathbb{R}}

\newcommand{\ZZ}{\mathbb{Z}}

\newcommand{\bA}{{\boldsymbol{A}}}

\newcommand{\ta}{\widetilde{a}}

\newcommand{\bB}{{\boldsymbol{B}}}
\newcommand{\tb}{\widetilde{b}}

\newcommand{\bC}{{\boldsymbol{C}}}

\newcommand{\bD}{{\boldsymbol{D}}}

\newcommand{\be}{{\boldsymbol{e}}}

\newcommand{\bh}{{\boldsymbol{h}}}
\newcommand{\bI}{{\boldsymbol{I}}}
\newcommand{\bJ}{{\boldsymbol{J}}}

\newcommand{\bM}{{\boldsymbol{M}}}

\newcommand{\bQ}{{\boldsymbol{Q}}}
\newcommand{\br}{{\boldsymbol{r}}}
\newcommand{\bS}{{\boldsymbol{S}}}

\newcommand{\bv}{{\boldsymbol{v}}}

\newcommand{\bx}{{\boldsymbol{x}}}

\newcommand{\bZ}{{\boldsymbol{Z}}}

\newcommand{\bDelta}{{\boldsymbol{\Delta}}}

\newcommand{\bfeta}{{\boldsymbol{\eta}}}
\newcommand{\btheta}{{\boldsymbol{\theta}}}
\newcommand{\tbtheta}{{\widetilde{\btheta}}}

\newcommand{\bttheta}{{\boldsymbol{\ttheta}}}

\newcommand{\bzero}{{\boldsymbol{0}}}

\newcommand{\cA}{{\mathcal A}}
\newcommand{\cB}{{\mathcal B}}

\newcommand{\cD}{{\mathcal D}}
\newcommand{\cE}{{\mathcal E}}

\newcommand{\cG}{{\mathcal G}}
\newcommand{\cF}{{\mathcal F}}

\newcommand{\cL}{{\mathcal L}}

\newcommand{\cN}{{\mathcal N}}

\newcommand{\cV}{{\mathcal V}}

\newcommand{\cX}{{\mathcal X}}
\newcommand{\cY}{{\mathcal Y}}
\newcommand{\cZ}{{\mathcal Z}}
\newcommand{\cW}{{\mathcal W}}

\newcommand{\tcY}{\widetilde{\cY}}

\newcommand{\tcV}{\widetilde{\cV}}

\newcommand{\dd}{\mathrm{d}}
\newcommand{\ee}{\mathrm{e}}

\newcommand{\ii}{\mathrm{i}}

\newcommand{\EE}{\operatorname{\mathbb{E}}}
\newcommand{\PP}{{\operatorname{\mathbb{P}}}}
\newcommand{\QQ}{{\operatorname{\mathbb{Q}}}}
\newcommand{\oo}{\operatorname{o}}
\newcommand{\OO}{\operatorname{O}}

\newcommand{\talpha}{\widetilde{\alpha}}
\newcommand{\tbeta}{\widetilde{\beta}}
\newcommand{\ttheta}{{\widetilde{\theta}}}

\newcommand{\vare}{\varepsilon}

\renewcommand{\mid}{\,|\,}

\renewcommand{\leq}{\leqslant}
\renewcommand{\geq}{\geqslant}

\newcommand{\stoch}{\stackrel{\PP}{\longrightarrow}}
\newcommand{\distr}{\stackrel{\cD}{\longrightarrow}}
\newcommand{\weak}{\Rightarrow}

\newcommand{\stochPT}{\stackrel{\PP_T}{\longrightarrow}}

\newcommand{\distre}{\stackrel{\cD}{=}}

\newcommand{\as}{\stackrel{{\mathrm{a.s.}}}{\longrightarrow}}

\newcommand{\bbone}{\mathbbm{1}}

\newcommand{\proofend}{\hfill\mbox{$\Box$}}

\numberwithin{equation}{section}

\theoremstyle{change} \theorembodyfont{\em}
\newtheorem{Lem}{Lemma.}[section]
\newtheorem{Thm}[Lem]{Theorem.}
\newtheorem{Pro}[Lem]{Proposition.}
\newtheorem{Cor}[Lem]{Corollary.}
\newtheorem{Def}[Lem]{Definition.}

\theorembodyfont{\rm}
\newtheorem{Rem}[Lem]{Remark.}

\begin{document}

\begin{center}
 {\bfseries\Large Local asymptotic quadraticity of statistical experiments} \\[2mm]
 {\bfseries\Large connected with a Heston model} \\[5mm]
 {\sc\large J\'anos Marcell Benke \ and \ Gyula $\text{Pap}^*$}
\end{center}

\vskip0.2cm

\noindent
 Bolyai Institute, University of Szeged,
 Aradi v\'ertan\'uk tere 1, H--6720 Szeged, Hungary.

\noindent e--mails: jbenke@math.u-szeged.hu (J. M. Benke),
                    papgy@math.u-szeged.hu (G. Pap).

\noindent * Corresponding author.

\renewcommand{\thefootnote}{}

\vspace*{-7mm}


\begin{abstract}
We study local asymptotic properties of likelihood ratios of certain
 Heston models.
We distinguish three cases: subcritical, critical and supercritical models. 
For the drift parameters, local asymptotic normality is proved in the
 subcritical case, only local asymptotic quadraticity is shown in the
 critical case, while in the supercritical case not even local asymptotic
 quadraticity holds.
For certain submodels, local asymptotic normality is proved in the
 critical case, and local asymptotic mixed normality is shown in the
 supercritical case.
As a consequence, asymptotically optimal (randomized) tests are constructed
 in cases of local asymptotic normality.
Moreover, local asymptotic minimax bound, and hence, asymptotic efficiency
 in the convolution theorem sense are concluded for the maximum likelihood
 estimators in cases of local asymptotic mixed normality.
\end{abstract}

\section{Introduction}

Heston models have been extensively used in financial mathematics since one 
 can well-fit them to real financial data set, and they are well-tractable
 from the point of view of computability as well, see Heston \cite{Hes}.

Let us consider a Heston model
 \begin{align}\label{Heston_SDE}
  \begin{cases}
   \dd Y_t = (a - b Y_t) \, \dd t + \sigma_1 \sqrt{Y_t} \, \dd W_t , \\
   \dd X_t = (\alpha - \beta Y_t) \, \dd t
             + \sigma_2  \sqrt{Y_t}
               \bigl(\varrho \, \dd W_t
                     + \sqrt{1 - \varrho^2} \, \dd B_t\bigr) ,
  \end{cases} \qquad t \geq 0 ,
 \end{align}
 where \ $a > 0$, \ $b, \alpha, \beta \in \RR$, \ $\sigma_1 > 0$,
 \ $\sigma_2 > 0$, \ $\varrho \in (-1, 1)$ \ and \ $(W_t, B_t)_{t\geq0}$ \ is a
 2-dimensional standard Wiener process.
Here one can interpret \ $X_t$ \ as the log-price of an asset, and \ $Y_t$ \ as
 the volatility of the asset price at time \ $t\geq0$.
\ The squared volatility process \ $(\sigma_2^2 Y_t)_{t\geq0}$ \ is a
 Cox--Ingersoll--Ross (CIR) process.
We distinguish three cases: subcritical if \ $b > 0$, \ critical if \ $b = 0$
 \ and supercritical if \ $b < 0$.
\ In this paper we study local asymptotic properties of the likelihood ratios
 of the model \eqref{Heston_SDE} concerning the drift parameter
 \ $(a, \alpha, b, \beta)$.

In case of the one-dimensional CIR process \ $Y$, \ Overbeck \cite{Ove} examined
 local asymptotic properties of the likelihood ratios concerning the drift
 parameter \ $(a, b)$, \ and proved the following results under the assumption
 \ $a \in \bigl(\frac{\sigma_1^2}{2}, \infty\bigr)$, \ which guarantees that the
 information matrix process tends to infinity almost surely.
It turned out that local asymptotic normality (LAN) is valid in the
 subcritical case.
In the critical case LAN has been proved for the submodel when \ $b = 0$ \ is
 known, and only local asymptotic quadraticity (LAQ) has been shown for the
 submodel when \ $a \in \bigl(\frac{\sigma_1^2}{2}, \infty\bigr)$ \ is known,
 but the asymptotic property of the experiment locally at \ $(a, 0)$ \ with a
 suitable two-dimensional localization sequence remained as an open question.
In the supercritical case local asymptotic mixed normality (LAMN) has been
 proved for the submodel when \ $a \in \bigl(\frac{\sigma_1^2}{2}, \infty\bigr)$
 \ is known.

For the Heston model \eqref{Heston_SDE}, we assume again
 \ $a \in \bigl(\frac{\sigma_1^2}{2}, \infty\bigr)$.
\ We prove LAN in the subcritical case (see Theorem \ref{Thm_subcrit}), LAQ in
 the critical case (see Theorem \ref{Thm_crit}), and show that LAQ does not hold
 in the supercritical case, although we can describe the asymptotic property of
 the experiment locally at \ $(a, \alpha, b, \beta)$ \ with a suitable
 four-dimensional degenerate localization sequence (see Theorem \ref{Thm_supercrit}). 
In the critical case LAN will be shown for the submodel when \ $b = 0$ \ and
 \ $\beta \in \RR$ \ are known (see Theorem \ref{Thm_crit}).
In the supercritical case LAMN will be proved for the submodel when
\ $a \in \bigl(\frac{\sigma_1^2}{2}, \infty\bigr)$ \ and \ $\alpha \in \RR$
 \ are known (see Theorem \ref{Thm_supercrit}).

If the LAN property holds then we obtain asymptotically optimal tests (see
 Remarks \ref{AOT_subcrit} and \ref{AOT_crit}) based on Theorem 15.4 and
 Addendum 15.5 of van der Vaart \cite{Vaart}.
 
If the LAMN property holds then we have a local asymptotic minimax bound for
 arbitrary estimators, see, e.g., Le Cam and Yang
 \cite[6.6, Theorem 1]{LeCamYang}.
Moreover, any maximum likelihood estimator attains this bound for bounded
 loss function (see Le Cam and Yang \cite[6.6, Remark 11]{LeCamYang}), and it
 is asymptotically efficient in H\'ajek's convolution theorem sense (for
 example, see, Le Cam and Yang \cite[6.6, Theorem 3 and Remark 13]{LeCamYang};
 Jeganathan \cite{Jeg}).
Asymptotic behavior of maximum likelihood estimators are described in all
 cases in Barczy and Pap \cite{BarPap}.

\section{Quadratic approximations to likelihood ratios}
\label{mlan}

Let \ $\NN$, \ $\ZZ_+$, \ $\RR$, \ $\RR_+$, \ $\RR_{++}$, \ $\RR_-$ \ and
 \ $\RR_{--}$ \ denote the sets of positive integers, non-negative integers,
 real numbers, non-negative real numbers, positive real numbers, non-positive
 real numbers and negative real numbers, respectively.
For \ $x , y \in \RR$, \ we will use the notations \ $x \land y := \min(x, y)$.
\ By \ $\|x\|$ \ and \ $\|A\|$, \ we denote the Euclidean norm of a vector
 \ $x \in \RR^d$ \ and the induced matrix norm of a matrix
 \ $A \in \RR^{d \times d}$, \ respectively.
By \ $\bI_d \in \RR^{d \times d}$, \ we denote the $d$-dimensional unit matrix.
In the sequel \ $\stoch$, \ $\distr$ \ and \ $\as$ \ will denote convergence
 in probability, in distribution and almost surely, respectively.

We recall some definitions and statements concerning quadratic approximations
 to likelihood ratios based on Jeganathan \cite{Jeg}, Le Cam and Yang
 \cite{LeCamYang} and van der Vaart \cite{Vaart}.

If \ $\PP$ \ and \ $\QQ$ \ are probability measures on a measurable space
 \ $(X, \cX)$, \ then
 \[
   \frac{\dd \PP}{\dd \QQ} : X \to \RR_+
 \]
 denotes the Radon--Nykodym derivative of the absolutely continuous part of
 \ $\PP$ \ with respect to \ $\QQ$.
\ If \ $(X, \cX, \PP)$ \ is a probability space and \ $(Y, \cY)$ \ is a
 measurable space, then the distribution of a measurable mapping
 \ $\xi : X \to Y$ \ under \ $\PP$ \ will be denoted by \ $\cL(\xi \mid \PP)$
 \ (i.e., \ $\cL(\xi \mid \PP)$ \ is the probability measure on \ $(Y, \cY)$
 \ defined by \ $\cL(\xi \mid \PP)(B) := \PP(\xi \in B)$, \ $B \in \cY$).

\begin{Def}
A statistical experiment is a triplet
 \ $\big( X, \cX, \{\PP_\theta : \theta \in \Theta\} \big)$, \ where
 \ $(X, \cX)$ \ is a measurable space and \ $\{\PP_\theta : \theta \in \Theta\}$
 \ is a family of probability measures on \ $(X, \cX)$. 
\ Its likelihood ratio process with base \ $\theta_0 \in \Theta$ \ is the
 stochastic process
 \[
   \left( \frac{\dd \PP_\theta}{\dd \PP_{\theta_0}} \right)_{\theta\in\Theta} .
 \]
\end{Def}

\begin{Def}
A family \ $(X_T, \cX_T, \{\PP_{\theta,T} : \theta \in \Theta\})_{T\in\RR_{++}}$
 \ of statistical experiments converges to a statistical experiment
 \ $(X, \cX, \{\PP_\theta : \theta \in \Theta\})$ \ as \ $T \to \infty$ \ if,
 for every finite subset \ $H \subset \Theta$ \ and every
 \ $\theta_0 \in \Theta$,
 \[
   \cL\left(\left( \frac{\dd \PP_{\theta,T}}{\dd \PP_{\theta_0,T}}
           \right)_{\theta\in H} \bigg| \, \PP_{\theta_0,T} \right)
   \weak
   \cL\left(\left( \frac{\dd\PP_\theta}{\dd\PP_{\theta_0}}
           \right)_{\theta\in H} \bigg| \, \PP_{\theta_0} \right)
  \qquad\text{as \ $T \to \infty$,}
 \]
 i.e., the finite dimensional distributions of the likelihood ratio process
 \ $\left( \frac{\dd \PP_{\theta,T}}{\dd \PP_{\theta_0,T}} \right)_{\theta\in\Theta}$
 \ under \ $\PP_{\theta_0,T}$ \ converges to the finite dimensional distributions
 of the likelihood ratio process
 \ $\left( \frac{\dd \PP_\theta}{\dd \PP_{\theta_0}} \right)_{\theta\in\Theta}$
 \ under \ $\PP_{\theta_0}$ \ as \ $T \to \infty$. 
\end{Def}

If \ $(X_T, \cX_T, \PP_T)$, \ $T \in \RR_{++}$, \ are probability spaces and
 \ $f_T : X_T \to \RR^p$, \ $T \in \RR_{++}$, \ are measurable functions, then
 \[
   f_T \stochPT 0 \quad \text{or} \quad f_T = \oo_{\PP_T}(1) \qquad
   \text{as \ $T \to \infty$}
 \]
 denotes convergence in \ $(\PP_T)_{T\in\RR_{++}}$-probabilities to 0 as
 \ $T \to \infty$, \ i.e., \ $\PP_T( \|f_T\| > \vare) \to 0$ \ as
 \ $T \to \infty$ \ for all \ $\vare \in \RR_{++}$.
\ Moreover,
 \[
   f_T = \OO_{\PP_T}(1) , \qquad T \in \RR_{++} ,
 \]
 denotes boundedness in \ $(\PP_T)_{T\in\RR_{++}}$-probabilities, \ i.e.,
 \ $\sup_{T\in\RR_{++}} \PP_T( \|f_T\| > K) \to 0$ \ as \ $K \to \infty$.

\begin{Rem}\label{oO}
Note that if \ $(\Omega, \cA, \PP)$ \ is a probability space and for each
 \ $T \in \RR_{++}$, \ $\xi_T : \Omega \to X_T$ \ is a random element with
 \ $\cL(\xi_T \mid \PP) = \PP_T$, \ then \ $f_T = \oo_{\PP_T}(1)$ \ as
 \ $T \to \infty$ \ or \ $f_T = \OO_{\PP_T}(1)$, \ $T \in \RR_{++}$, \ if and
 only if \ $f_T \circ \xi_T = \oo_\PP(1)$ \ as \ $T \to \infty$ \ or
 \ $f_T \circ \xi_T = \OO_\PP(1)$, \ $T \in \RR_{++}$, \ respectively.
Indeed, \ $\PP_T( \|f_T\| > c) = \PP( \|f_T(\xi_T)\| > c)$ \ for all
 \ $T \in \RR_{++}$ \ and all \ $c \in \RR_{++}$.
\ Moreover, \ $f_T = \OO_{\PP_T}(1)$, \ $T \in \RR_{++}$, \ if and only if the
 family \ $(\cL(f_T \mid \PP_T))_{T\in\RR_{++}}$ \ of probability measures is
 tight, and hence, for each sequence \ $T_n \in \RR_{++}$, \ $n \in \NN$,
 \ with \ $T_n \to \infty$ \ as \ $n \to \infty$, \ there exist a subsequence
 \ $T_{n_k}$, \ $k \in \NN$, \ and a probability measure \ $\mu$ \ on
 \ $(\RR^p, \cB(\RR^p))$, \ such that
 \ $\cL(f_{T_{n_k}} \mid \PP_{T_{n_k}}) \weak \mu$ \ as \ $k \to \infty$.
\ In this case, \ $\mu$ \ is called an accumulation point of the family
 \ $(\cL(f_T \mid \PP_T))_{T\in\RR_{++}}$.
\proofend
\end{Rem}

\begin{Def}
Let \ $\Theta \subset \RR^p$ \ be an open set.
A family \ $(X_T, \cX_T, \{\PP_{\btheta,T} : \btheta \in \Theta\})_{T\in\RR_{++}}$
 \ of statistical experiments is said to have locally asymptotically quadratic
 (LAQ) likelihood ratios at \ $\btheta \in \Theta$ \ if there exist (scaling)
 matrices \ $\br_{\btheta,T} \in \RR^{p\times p}$, \ $T \in \RR_{++}$, \ measurable
 functions (statistics) \ $\bDelta_{\btheta,T} : X_T \to \RR^p$,
 \ $T \in \RR_{++}$, \ and \ $\bJ_{\btheta,T} : X_T \to \RR^{p \times p}$,
 \ $T \in \RR_{++}$, \ such that
 \begin{equation}\label{LAQ}
   \log \frac{\dd\PP_{\btheta+\br_{\btheta,T}\bh_T,T}}{\dd\PP_{\btheta,T}}
   = \bh_T^\top \bDelta_{\btheta,T}- \frac{1}{2} \bh_T^\top \bJ_{\btheta,T} \bh_T
     + \oo_{\PP_{\btheta,T}}(1)
   \qquad\text{as \ $T \to \infty$}
 \end{equation}
 whenever \ $\bh_T \in \RR^p$, \ $T \in \RR_{++}$, \ is a bounded family
 satisfying \ $\btheta + \br_{\btheta,T} \bh_T \in \Theta$ \ for all
 \ $T \in \RR_{++}$,
 \begin{equation}\label{LAQO}
  (\bDelta_{\btheta,T}, \bJ_{\btheta,T}) = \OO_{\PP_{\btheta,T}}(1) ,
  \qquad T \in \RR_{++} ,
 \end{equation}
 and for each accumulation point \ $\mu_\btheta$ \ of the family
 \ $(\cL( (\bDelta_{\btheta,T}, \bJ_{\btheta,T}) \mid \PP_{\btheta,T}))_{T\in\RR_{++}}$
 \ as \ $T \to \infty$, \ which is a probability measure on
 \ $(\RR^p \times \RR^{p \times p}, \cB(\RR^p \times \RR^{p \times p}))$, \ we have
 \begin{equation}\label{LAQJ}
  \mu_\btheta\left(\left\{(\bDelta, \bJ) \in \RR^p \times \RR^{p \times p}
                         : \text{$\bJ$ is symmetric and strictly positive
                                 definite}
                  \right\}\right)
  = 1   
 \end{equation}
 and
 \begin{equation}\label{LAQDJ}
  \int_{\RR^p \times \RR^{p \times p}}
   \exp\left\{ \bh^\top \bDelta - \frac{1}{2} \bh^\top \bJ \bh \right\}
   \, \mu_\btheta(\dd\bDelta, \dd\bJ)
  = 1 
 \end{equation}
 whenever \ $\bh \in \RR^p$ \ such that there exist \ $T_k \in \RR_{++}$,
 \ $k \in \NN$, \ and \ $\bh_{T_k} \in \RR^p$,
 \ $k \in \NN$, \ with \ $\bh_{T_k} \to \bh$ \ as \ $k \to \infty$,
 \ $\btheta + \br_{\btheta,T_k} \bh_{T_k} \in \Theta$ \ for all \ $k \in \NN$.
\end{Def}

\begin{Def}
Let \ $\Theta \subset \RR^p$ \ be an open set.
A family \ $(X_T, \cX_T, \{\PP_{\btheta,T} : \btheta \in \Theta\})_{T\in\RR_{++}}$
 \ of statistical experiments is said to have locally asymptotically mixed
 normal (LAMN) likelihood ratios at \ $\btheta \in \Theta$ \ if it is LAQ at
 \ $\btheta \in \Theta$, \ and for each accumulation point \ $\mu_\btheta$ \ of
 the family
 \ $(\cL( (\bDelta_{\btheta,T}, \bJ_{\btheta,T}) \mid \PP_{\btheta,T}))_{T\in\RR_{++}}$
 \ as \ $T \to \infty$, \ we have
 \[
   \int_{\RR^p \times B} \ee^{\ii\,\bh^\top\!\bDelta} \, \mu_\btheta(\dd\bDelta, \dd\bJ)
   = \int_{\RR^p \times B} \ee^{-\bh^\top\!\bJ\bh/2} \, \mu_\btheta(\dd\bDelta, \dd\bJ) ,
   \qquad B \in \cB(\RR^{p \times p}), \quad \bh \in \RR^p ,
 \]
 i.e., the conditional distribution of \ $\bDelta$ \ given \ $\bJ$ \ under
 \ $\mu_\btheta$ \ is \ $\cN_p(\bzero, \bJ)$, \ or, equivalently, 
 \ $\mu_\btheta = \cL((\eta_\btheta \cZ, \eta_\btheta \eta_\btheta^\top) \mid \PP)$,
 \ where \ $\cZ : \Omega \to \RR^p$ \ and
 \ $\eta_\btheta : \Omega \to \RR^{p \times p}$ \ are independent random elements
 on a probability space \ $(\Omega, \cF, \PP)$ \ such that
 \ $\cL(\cZ \mid \PP) = \cN_p(\bzero, \bI_p)$.
\end{Def}

\begin{Def}
Let \ $\Theta \subset \RR^p$ \ be an open set.
A family \ $(X_T, \cX_T, \{\PP_{\btheta,T} : \btheta \in \Theta\})_{T\in\RR_{++}}$
 \ of statistical experiments is said to have locally asymptotically normal
 (LAN) likelihood ratios at \ $\btheta \in \Theta$ \ if it is LAMN at
 \ $\btheta \in \Theta$, \ and for each accumulation point \ $\mu_\btheta$ \ of
 the family
 \ $(\cL( (\bDelta_{\btheta,T}, \bJ_{\btheta,T}) \mid \PP_{\btheta,T}))_{T\in\RR_{++}}$
 \ as \ $T \to \infty$, \ we have
 \[
   \mu_\btheta = \cN_p(\bzero, \bJ_\btheta) \times \delta_{\bJ_\btheta}
 \]
 with some symmetric, strictly positive definite matrix
 \ $\bJ_\btheta \in \RR^{p \times p}$, \ where \ $\delta_{\bJ_\btheta}$ \ denotes the
 Dirac measure on \ $(\RR^{p \times p}, \cB(\RR^{p \times p}))$, \ concentrated in
 \ $\bJ_\btheta$.
\end{Def}

We will need Le Cam's first lemma, see, e.g, Lemma 6.4 in van der Vaart
 \cite{Vaart}.
We start with the definition of contiguity of families of probability measures.

\begin{Def}
Let \ $(X_T, \cX_T)$, \ $T \in \RR_{++}$, \ be measurable spaces.
For each \ $T \in \RR_{++}$, \ let \ $\PP_T$ \ and \ $\QQ_T$ \ be probability
 measures on \ $(X_T, \cX_T)$.
\ The family \ $(\QQ_T)_{T\in\RR_{++}}$ \ is said to be contiguous with respect to
 the family \ $(\PP_T)_{T\in\RR_{++}}$ \ if \ $\QQ_T(A_T) \to 0$ \ as
 \ $T \to \infty$ \ whenever \ $A_T \in \cX_T$, \ $T \in \RR_{++}$, \ such that
 \ $\PP_T(A_T) \to 0$ \ as \ $T \to \infty$.
\ This will be denoted by
 \ $(\QQ_T)_{T\in\RR_{++}} \triangleleft (\PP_T)_{T\in\RR_{++}}$.
\ The families \ $(\PP_T)_{T\in\RR_{++}}$ \ and \ $(\QQ_T)_{T\in\RR_{++}}$ \ are said
 to be mutually contiguous if both
 \ $(\PP_T)_{T\in\RR_{++}} \triangleleft (\QQ_T)_{T\in\RR_{++}}$ \ and
 \ $(\QQ_T)_{T\in\RR_{++}} \triangleleft (\PP_T)_{T\in\RR_{++}}$ \ hold.
\end{Def}

\begin{Lem}[Le Cam's first lemma]\label{LeCam1}
Let \ $(X_T, \cX_T)$, \ $T \in \RR_{++}$, \ be measurable spaces.
For each \ $T \in \RR_{++}$, \ let \ $\PP_T$ \ and \ $\QQ_T$ \ be probability
 measures on \ $(X_T, \cX_T)$.
\ Then the following statements are equivalent:
 \renewcommand{\labelenumi}{{\rm(\roman{enumi})}}
 \begin{enumerate}
  \item
   $(\QQ_T)_{T\in\RR_{++}} \triangleleft (\PP_T)_{T\in\RR_{++}}$;
  \item
   If
    \ $\cL\Bigl( \frac{\dd\PP_{T_k}}{\dd\QQ_{T_k}} \, \Big| \, \QQ_{T_k} \Bigr)
       \weak \nu$
    \ as \ $k \to \infty$ \ for some sequence \ $(T_k)_{k\in\NN}$ \ with
    \ $T_k \to \infty$ \ as \ $T \to \infty$, \ where \ $\nu$ \ is a
    probability measure on \ $(\RR_+, \cB(\RR_+))$, \ then
    \ $\nu(\RR_{++}) = 1$;
  \item
   If
    \ $\cL\Bigl( \frac{\dd\QQ_{T_k}}{\dd\PP_{T_k}} \, \Big| \, \PP_{T_k} \Bigr)
       \weak \mu$
    \ as \ $k \to \infty$ \ for some sequence \ $(T_k)_{k\in\NN}$ \ with
    \ $T_k \to \infty$ \ as \ $T \to \infty$, \ where \ $\mu$ \ is a
    probability measure on \ $(\RR_+, \cB(\RR_+))$, \ then
    \ $\int_{\RR_{++}} x \, \mu(\dd x) = 1$;
  \item
   $\cL(f_T \mid \QQ_T) \weak 0$ \ as \ $T \to \infty$ \ whenever
    \ $f_T : X_T \to \RR^p$, \ $T \in \RR_{++}$, \ are measurable functions and
    \ $\cL(f_T \mid \PP_T) \weak 0$ \ as \ $T \to \infty$.
 \end{enumerate}
\end{Lem}

We will need a version of general form of Le Cam's third lemma, which is
 Theorem 6.6 in van der Vaart \cite{Vaart}.

\begin{Thm}\label{LeCam3}
Let \ $(X_T, \cX_T)$, \ $T \in \RR_{++}$, \ be measurable spaces.
For each \ $T \in \RR_{++}$, \ let \ $\PP_T$ \ and \ $\QQ_T$ \ be probability
 measures on \ $(X_T, \cX_T)$.
\ Let \ $f_T : X_T \to \RR^p$, \ $T \in \RR_{++}$, \ be measurable functions.
Suppose that the family \ $(\QQ_T)_{T\in\RR_{++}}$ \ is contiguous with respect to
 the family \ $(\PP_T)_{T\in\RR_{++}}$ \ and
 \[
   \cL\left( \left(f_T, \frac{\dd\QQ_T}{\dd\PP_T}\right)
             \, \bigg| \, \PP_T \right)
   \weak \nu \qquad \text{as \ $T \to \infty$,}
 \]
 where \ $\nu$ \ is a probability measure on
 \ $(\RR^p \times \RR_+, \cB(\RR^p \times \RR_+))$.
\ Then \ $\cL(f_T \mid \QQ_T) \weak \mu$ \ as \ $T \to \infty$, \ where \ $\mu$
 \ is the probability measure on \ $(\RR^p, \cB(\RR^p))$ \ given by
 \[
   \mu(B) := \int_{\RR^p \times \RR_+} \bbone_B(f) V \, \nu(\dd f, \dd V) , \qquad
   B \in \cB(\RR^p) .
 \]
\end{Thm}

The following convergence theorem is Proposition 1 in Jeghanathan \cite{Jeg}.
In fact, it is a generalization of Theorems 9.4 and 9.8 of van der Vaart
 \cite{Vaart}, which are valid for LAMN and LAN families of experiments.
For completeness, we give a proof.

\begin{Thm}\label{conv_LAQ}
Let \ $\Theta \subset \RR^p$ \ be an open set.
Let \ $(X_T, \cX_T, \{\PP_{\btheta,T} : \btheta \in \Theta\})_{T\in\RR_{++}}$ \ be a
 family of statistical experiments.
Assume that LAQ is satisfied at \ $\btheta \in \Theta$.
\ Let \ $T_k \in \RR_{++}$, \ $k \in \NN$, \ be such that \ $T_k \to \infty$
 \ and
 \ $\cL( (\bDelta_{\btheta,T_k}, \bJ_{\btheta,T_k}) \mid \PP_{\btheta,T_k})
    \weak \mu_\btheta$
 \ as \ $k \to \infty$.
\ Then, for every \ $\bh_{T_k} \in \RR^p$, \ $k \in \NN$, \ with
 \ $\bh_{T_k} \to \bh$\ as \ $k \to \infty$ \ and
 \ $\btheta + \br_{\btheta,T_k} \bh_{T_k} \in \Theta$ \ for all \ $k \in \NN$,
 \ we have
 \ $\cL( (\bDelta_{\btheta,T_k}, \bJ_{\btheta,T_k})
         \mid \PP_{\btheta+\br_{\btheta,T_k}\bh_{T_k},T_k})
    \weak \QQ_{\btheta,\bh}$
 \ as \ $k \to \infty$, \ where
 \[
   \QQ_{\btheta,\bh}(B)
   := \int_B
       \exp\left\{ \bh^\top \bDelta - \frac{1}{2} \bh^\top \bJ \bh \right\}
       \, \mu_\btheta(\dd\bDelta, \dd\bJ) , \qquad
   B \in \cB(\RR^p \times \RR^{p\times p}) .
 \]
Consequently, the sequence
 \ $(X_{T_k}, \cX_{T_k},
     \{\PP_{\btheta+\br_{\btheta,T_k}\bh,T_k} : \bh \in \RR^p\})_{k\in\NN}$
 \ of statistical experiments converges to the statistical experiment
 \ $(\RR^p \times \RR^{p\times p}, \cB(\RR^p \times \RR^{p\times p}),
     \{\QQ_{\btheta,\bh} : \bh \in \RR^p\})$
 \ as \ $k \to \infty$.
\end{Thm}

Note that for each \ $\bh \in \RR^p$, \ the probability measures
 \ $\QQ_{\btheta,\bh}$ \ and \ $\QQ_{\btheta,\bzero}$ \ are equivalent, and
 \[
   \frac{\dd\QQ_{\btheta,\bh}}{\dd\QQ_{\btheta,\bzero}}(\bDelta, \bJ)
   = \exp\left\{ \bh^\top \bDelta - \frac{1}{2} \bh^\top \bJ \bh \right\} ,
   \qquad (\bDelta, \bJ) \in \RR^p \times \RR^{p\times p} .
 \]

\noindent{\bf Proof.}
Let \ $(\Omega, \cA, \PP)$ \ be a probability space and let
 \ $(\bDelta, \bJ) : \Omega \to \RR^p \times \RR^{p\times p}$ \ be a measurable
 function such that \ $\cL( (\bDelta, \bJ) \mid \PP) = \mu_\btheta$.
\ Using
 \ $\cL( (\bDelta_{\btheta,T_k}, \bJ_{\btheta,T_k}) \mid \PP_{\btheta,T_k})
    \weak \mu_\btheta$
 \ as \ $k \to \infty$, \ by Slutsky's lemma,
 \[
   \cL\left( \frac{\dd\PP_{\btheta+\br_{\btheta,T_k}\bh,T_k}}{\dd\PP_{\btheta,T_k}}
             \, \bigg| \, \PP_{\btheta,T_k} \right)
   \weak \cL\left(\exp\left\{ \bh^\top \bDelta
                              - \frac{1}{2} \bh^\top \bJ \bh \right\}\right)
   \qquad \text{as \ $T \to \infty$.}
 \]
By \eqref{LAQJ} and \eqref{LAQDJ}, applying Lemma \ref{LeCam1}, we conclude
 that the sequences \ $(\PP_{\btheta+\br_{\btheta,T_k}\bh,T_k})_{k\in\NN}$ \ and
 \ $(\PP_{\btheta,T_k})_{k\in\NN}$ \ are mutually contiguous.
Therefore, for each \ $\bh, \bh_0 \in \RR^p$, \ the probability of the set on
 which we have
 \[
   \log \frac{\dd\PP_{\btheta+\br_{\btheta,T_k}\bh,T_k}}
             {\dd\PP_{\btheta+\br_{\btheta,T_k}\bh_0,T_k}}
   = \log \frac{\dd\PP_{\btheta+\br_{\btheta,T_k}\bh,T_k}}{\dd\PP_{\btheta,T_k}}
     - \log \frac{\dd\PP_{\btheta+\br_{\btheta,T_k}\bh_0,T_k}}{\dd\PP_{\btheta,T_k}} ,
 \]
 converges to one.
By \eqref{LAQ}, we obtain
 \[
   \log \frac{\dd\PP_{\btheta+\br_{\btheta,T_k}\bh,T_k}}
             {\dd\PP_{\btheta+\br_{\btheta,T_k}\bh_0,T_k}}
   = (\bh - \bh_0)^\top \bDelta_{\btheta,T_k}
     - \frac{1}{2} \bh^\top \bJ_{\btheta,T_k} \bh
     + \frac{1}{2} \bh_0^\top \bJ_{\btheta,T_k} \bh_0
     + \oo_{\PP_{\btheta,T}}(1)
   \qquad\text{as \ $T \to \infty$.}
 \]
Hence it suffices to observe that
 \ $\cL( (\bDelta_{\btheta,T_k}, \bJ_{\btheta,T_k})
         \mid \PP_{\btheta+\br_{\btheta,T_k}\bh,T_k})
    \weak \QQ_{\btheta,\bh}$
 \ as \ $k \to \infty$ \ for all \ $\bh \in \RR^p$ \ follows from Theorem
 \ref{LeCam3}.
\proofend

The following statements are trivial consequences of Theorem \ref{conv_LAQ},
 and they can also be derived from Theorems 9.4 and 9.8 of van der Vaart
 \cite{Vaart}.

\begin{Pro}\label{conv_LAMN}
Let \ $\Theta \subset \RR^p$ \ be an open set.
Let \ $(X_T, \cX_T, \{\PP_{\btheta,T} : \btheta \in \Theta\})_{T\in\RR_{++}}$ \ be a
 family of statistical experiments.
Assume that LAMN is satisfied at \ $\btheta \in \Theta$.
\ Let \ $T_k \in \RR_{++}$, \ $k \in \NN$, \ be such that
 \ $\cL( (\bDelta_{\btheta,T_k}, \bJ_{\btheta,T_k}) \mid \PP_{\btheta,T_k})
    \weak \cL((\eta_\btheta \cZ, \eta_\btheta \eta_\btheta^\top) \mid \PP)$
 \ as \ $k \to \infty$, \ where \ $\cZ : \Omega \to \RR^p$ \ and
 \ $\eta_\btheta : \Omega \to \RR^{p \times p}$ \ are independent random elements
 on a probability space \ $(\Omega, \cF, \PP)$ \ such that
 \ $\cL(\cZ \mid \PP) = \cN_p(\bzero, \bI_p)$.
\ Then, for every \ $\bh_{T_k} \in \RR^p$, \ $k \in \NN$, \ with
 \ $\bh_{T_k} \to \bh$\ as \ $k \to \infty$ \ and
 \ $\btheta + \br_{\btheta,T_k} \bh_{T_k} \in \Theta$ \ for all \ $k \in \NN$,
 \ we have
 \ $\cL( (\bDelta_{\btheta,T_k}, \bJ_{\btheta,T_k})
         \mid \PP_{\btheta+\br_{\btheta,T_k}\bh_{T_k},T_k})
    \weak \cL((\eta_\btheta \cZ + \eta_\btheta \eta_\btheta^\top \bh,
               \eta_\btheta \eta_\btheta^\top)
              \mid \PP)$
 \ as \ $k \to \infty$.
\ Consequently, the sequence
 \ $(X_{T_k}, \cX_{T_k},
     \{\PP_{\btheta+\br_{\btheta,T_k}\bh,T_k} : \bh \in \RR^p\})_{T\in\RR_{++}}$
 \ of statistical experiments converges to the statistical experiment
 \ $(\RR^p \times \RR^{p\times p}, \cB(\RR^p \times \RR^{p\times p}),
     \{\cL((\eta_\btheta \cZ + \eta_\btheta \eta_\btheta^\top \bh,
            \eta_\btheta \eta_\btheta^\top)
           \mid \PP) : \bh \in \RR^p\})$
 \ as \ $k \to \infty$.
\end{Pro}

\begin{Pro}\label{conv_LAN}
Let \ $\Theta \subset \RR^p$ \ be an open set.
Let \ $(X_T, \cX_T, \{\PP_{\btheta,T} : \btheta \in \Theta\})_{T\in\RR_{++}}$ \ be a
 family of statistical experiments.
Assume that LAN is satisfied at \ $\btheta \in \Theta$.
\ Let \ $T_k \in \RR_{++}$, \ $k \in \NN$, \ be such that
 \ $\cL( (\bDelta_{\btheta,T_k}, \bJ_{\btheta,T_k}) \mid \PP_{\btheta,T_k})
    \weak \cN_p(\bzero, \bJ_\btheta) \times \delta_{\bJ_\btheta}$
 \ as \ $k \to \infty$ \ with some symmetric, strictly positive definite
 matrix \ $\bJ_\btheta \in \RR^{p \times p}$.
\ Then, for every \ $\bh_{T_k} \in \RR^p$, \ $k \in \NN$, \ with
 \ $\bh_{T_k} \to \bh$\ as \ $k \to \infty$ \ and
 \ $\btheta + \br_{\btheta,T_k} \bh_{T_k} \in \Theta$ \ for all \ $k \in \NN$,
 \ we have
 \ $\cL( (\bDelta_{\btheta,T_k}, \bJ_{\btheta,T_k})
         \mid \PP_{\btheta+\br_{\btheta,T_k}\bh_{T_k},T_k})
    \weak \cN_p(\bJ_\btheta \bh, \bJ_\btheta) \times \delta_{\bJ_\btheta}$
 \ as \ $k \to \infty$.
\ Consequently, the sequence
 \ $(X_{T_k}, \cX_{T_k},
     \{\PP_{\btheta+\br_{\btheta,T_k}\bh,T_k} : \bh \in \RR^p\})_{T\in\RR_{++}}$
 \ of statistical experiments converges to the statistical experiment
 \ $(\RR^p, \cB(\RR^p), \{\cN_p(\bJ_\btheta \bh, \bJ_\btheta) : \bh \in \RR^p\})$
 \ as \ $k \to \infty$.
\end{Pro}

\section{Asymptotically optimal tests}
\label{AOT}

\begin{Def}
A (randomized) test (function) in a statistical experiment 
 \ $(X, \cX, \{\PP_\btheta : \btheta \in \Theta\})$ \ is a Borel measurable
 function \ $\phi : X \to [0,1]$.
\ (The interpretation is that if \ $x \in X$ \ is observed, then a null
 hypothesis \ $H_0 \subset \Theta$ \ is rejected with probability
 \ $\phi(x)$.)

The power function of a test \ $\phi$ \ is the function
 \ $\theta \mapsto \int_X \phi(x) \, \PP_\btheta(\dd x)$. 
\ (This gives the probability that the null hypothesis \ $H_0$ \ is rejected.)

For \ $\alpha \in (0, 1)$, \ a test \ $\phi$ \ is of level \ $\alpha$ \ for
 testing a null hypothesis \ $H_0$ \ if
 \[
   \sup\left\{\int_X \phi(x) \, \PP_\btheta(\dd x) : \theta \in H_0\right\}
   \leq \alpha .
 \]
\end{Def}

If the LAN property holds then one obtains asymptotically optimal tests in the
 following way, see, e.g., Theorem 15.4 and Addendum 15.5 of van der Vaart \cite{Vaart}.

\begin{Thm}\label{AOT_LAN}
Let \ $\Theta \subset \RR^p$ \ be an open set.
Let \ $(X_T, \cX_T, \{\PP_{\btheta,T} : \btheta \in \Theta\})_{T\in\RR_{++}}$ \ be a
 family of statistical experiments such that LAN is satisfied at
 \ $\btheta_0 \in \Theta$.
\ Let \ $T_k \in \RR_{++}$, \ $k \in \NN$, \ be such that
 \ $\cL( (\bDelta_{\btheta_0,T_k}, \bJ_{\btheta_0,T_k}) \mid \PP_{\btheta_0,T_k})
    \weak \cN_p(\bzero, \bJ_{\btheta_0}) \times \delta_{\bJ_{\btheta_0}}$
 \ as \ $k \to \infty$ \ with some symmetric, strictly positive definite
 matrix \ $\bJ_{\btheta_0} \in \RR^{p \times p}$.
\ Let \ $\psi : \Theta \to \RR$ \ be differentiable at
 \ $\btheta_0 \in \Theta$ \ with \ $\psi(\btheta_0) = 0$ \ and
 \ $\psi'(\btheta_0) \ne \bzero$. 
\ Let \ $\alpha \in (0, 1)$.
\ For each \ $k \in \NN$, \ let \ $\phi_k : X_{T_k} \to [0,1]$ \ be a test of
 level \ $\alpha$ \ for testing \ $H_0 : \psi(\btheta) \leq 0$ \ against
 \ $H_1 : \psi(\btheta) > 0$, \ i.e., it is a Borel measurable function such
 that
 \[
   \sup\left\{\int_{X_{T_k}} \phi_k(x) \, \PP_{\btheta,T_k}(\dd x)
              : \btheta \in \Theta , \, \psi(\btheta) \leq 0\right\}
   \leq \alpha .
 \]
Then for each \ $\bh \in \RR^p$ \ with
 \ $\langle \psi'(\btheta_0), \bh \rangle > 0$, \ the power function of the
 test \ $\phi_k$ \ satisfies
 \[
   \limsup_{k\to\infty}
    \int_{X_{T_k}} \phi_k(x) \, \PP_{\btheta_0+\br_{\btheta_0,T_k}\bh,T_k}(\dd x)
   \leq 1
        - \Phi\left(z_\alpha
                    - \frac{\langle \psi'(\btheta_0), \bh \rangle}
                           {\sqrt{\langle \bJ_{\btheta_0}^{-1} \psi'(\btheta_0),
                                          \psi'(\btheta_0) \rangle}}\right) ,
 \]
 where \ $\Phi$ \ denotes the standard normal distribution function, and
 \ $z_\alpha$ \ denotes the upper \ $\alpha$-quantile of the standard normal
 distribution.
 
Moreover, if \ $S_{\btheta_0,k} : X_{T_k} \to \RR$, \ $k \in \NN$, \ are Borel
 measurable functions such that
 \[
   S_{\btheta_0,k}
   = \frac{\langle \bJ_{\btheta_0}^{-1} \bDelta_{\btheta_0,T_k},
                   \psi'(\btheta_0) \rangle}
          {\sqrt{\langle \bJ_{\btheta_0}^{-1} \psi'(\btheta_0),
                                          \psi'(\btheta_0) \rangle}}
     + \oo_{\PP_{\btheta_0,k}}(1) , \qquad k \in \NN ,
 \]
 then the family of tests that reject for values \ $S_{\btheta_0,k}$ \ exceeding
 \ $z_\alpha$ \ is asymptotically optimal for testing
 \ $H_0 : \psi(\btheta) \leq 0$ \ against \ $H_1 : \psi(\btheta) > 0$ \ in the
 sense that for every \ $\bh \in \RR^p$ \ with
 \ $\langle \psi'(\btheta_0), \bh \rangle > 0$,
 \[
  \PP\big(S_{\btheta_0,k}(x) \geq z_\alpha\big)
  \to 1
      - \Phi\left(z_\alpha
                  - \frac{\langle \psi'(\btheta_0, \bh \rangle}
                         {\sqrt{\langle \bJ_{\btheta_0}^{-1} \psi'(\btheta_0),
                                        \psi'(\btheta_0) \rangle}}\right)
  \qquad\text{as \ $k \to \infty$.}
 \]
\end{Thm}

\section{Local asymptotic minimax bound for estimators}
\label{LAMB}

If LAMN property holds then we have the following local asymptotic minimax bound for
 arbitrary estimators, see, e.g., Le Cam and Yang \cite[6.6, Theorem 1]{LeCamYang}.

\begin{Pro}\label{LAMB_LAMN}
Let \ $\Theta \subset \RR^p$ \ be an open set.
Let \ $(X_T, \cX_T, \{\PP_{\btheta,T} : \btheta \in \Theta\})_{T\in\RR_{++}}$ \ be a
 family of statistical experiments.
Assume that LAMN is satisfied at \ $\btheta \in \Theta$.
\ Let \ $T_k \in \RR_{++}$, \ $k \in \NN$, \ be such that
 \ $\cL( (\bDelta_{\btheta,T_k}, \bJ_{\btheta,T_k}) \mid \PP_{\btheta,T_k})
    \weak \cL((\eta_\btheta \cZ, \eta_\btheta \eta_\btheta^\top) \mid \PP)$
 \ as \ $k \to \infty$, \ where \ $\cZ : \Omega \to \RR^p$ \ and
 \ $\eta_\btheta : \Omega \to \RR^{p \times p}$ \ are independent random elements
 on a probability space \ $(\Omega, \cF, \PP)$ \ such that
 \ $\cL(\cZ \mid \PP) = \cN_p(\bzero, \bI_p)$.
\ Let \ $w : \RR^p \to \RR_+$ \ be a bowl-shaped loss function, i.e., for each
 \ $c \in \RR_+$, \ the set \ $\{x \in \RR^p : w(x) \leq c\}$ \ is closed, convex
 and symmetric.
Then, for arbitrary estimators (statistics, i.e., measurable functions)
 \ $\tbtheta_T : X_T \to \RR^p$, \ $T \in \RR_+$, \ of the parameter \ $\btheta$,
 \ one has
 \[
   \lim_{c\to\infty}
    \liminf_{k\to\infty}
     \sup_{\{x\in X_{T_k}:\|r_{\btheta,T_k}^{-1}(\tbtheta_{T_k}(x)-\btheta)\|\leq c\}}
      \int_{X_{T_k}}
       w\bigl(r_{\btheta,T_k}^{-1}(\tbtheta_{T_k}(x)-\btheta)\bigr)
       \, \PP_{\btheta,T_k}(\dd x)
   \geq
   \EE\bigl[w((\eta_\btheta^\top)^{-1}\cZ)\bigr] .
 \]
\end{Pro}

Maximum likelihood estimators attain this bound for bounded loss function \ $w$,
 \ see, e.g., Le Cam and Yang \cite[6.6, Remark 11]{LeCamYang}. 
Moreover, maximum likelihood estimators are asymptotically efficient in H\'ajek's
 convolution theorem sense (for example, see, Le Cam and Yang
 \cite[6.6, Theorem 3 and Remark 13]{LeCamYang}; Jeganathan \cite{Jeg}).

\section{Heston models}
\label{Prel}

The next proposition is about the existence and uniqueness of a strong
 solution of the SDE \eqref{Heston_SDE}, see, e.g., Barczy and Pap
 \cite[Proposition 2.1]{BarPap}.

\begin{Pro}\label{Pro_Heston}
Let \ $\bigl(\Omega, \cF, \PP\bigr)$ \ be a probability space.
Let \ $(W_t, B_t)_{t\in\RR_+}$ \ be a 2-dimensional standard Wiener process.
Let \ $(\eta_0, \zeta_0)$ \ be a random vector independent of
 \ $(W_t, B_t)_{t\in\RR_+}$ \ satisfying \ $\PP(\eta_0 \in \RR_+) = 1$.
\ Then for all \ $a \in \RR_{++}$, \ $b, \alpha, \beta \in \RR$,
 \ $\sigma_1, \sigma_2 \in \RR_{++}$, \ $\varrho \in (-1, 1)$, \ there is a
 (pathwise) unique strong solution \ $(Y_t, X_t)_{t\in\RR_+}$ \ of the SDE
 \eqref{Heston_SDE} such that \ $\PP((Y_0, X_0) = (\eta_0, \zeta_0)) = 1$ \ and
 \ $\PP(\text{$Y_t \in \RR_+$ \ for all \ $t \in \RR_+$}) = 1$.
\end{Pro}

Based on the asymptotic behavior of the expectations \ $(\EE(Y_t), \EE(X_t))$
 \ as \ $t \to \infty$, \ one can classify Heston processes given by the SDE
 \eqref{Heston_SDE}, see Barczy and Pap \cite{BarPap}.

\begin{Def}\label{Def_criticality}
Let \ $(Y_t, X_t)_{t\in\RR_+}$ \ be the unique strong solution of the SDE
 \eqref{Heston_SDE} satisfying \ $\PP(Y_0 \in \RR_+) = 1$.
\ We call \ $(Y_t, X_t)_{t\in\RR_+}$ \ subcritical, critical or supercritical if
 \ $b \in \RR_{++}$, \ $b = 0$ \ or \ $b \in \RR_{--}$, \ respectively.
\end{Def}

The following result states ergodicity of the process \ $(Y_t)_{t\in\RR_+}$
 \ given by the first equation in \eqref{Heston_SDE} in the subcritical case,
 see, e.g., Cox et al.\ \cite[Equation (20)]{CoxIngRos},
 Li and Ma \cite[Theorem 2.6]{LiMa} or Theorem 4.1 in Barczy et al.\
 \cite{BarDorLiPap2}.

\begin{Thm}\label{Ergodicity}
Let \ $a, b, \sigma_1 \in \RR_{++}$.
\ Let \ $(Y_t)_{t\in\RR_+}$ \ be the unique strong solution of the first equation
 of the SDE \eqref{Heston_SDE} satisfying \ $\PP(Y_0 \in \RR_+) = 1$.
 \renewcommand{\labelenumi}{{\rm(\roman{enumi})}}
 \begin{enumerate}
  \item
   Then \ $Y_t \distr Y_\infty$ \ as \ $t \to \infty$, \ and the distribution of
    \ $Y_\infty$ \ is given by
    \begin{align}\label{Laplace}
     \EE(\ee^{-\lambda Y_\infty})
     = \biggl(1 + \frac{\sigma_1^2}{2b} \lambda\biggr)^{-2a/\sigma_1^2} ,
     \qquad \lambda \in \RR_+ ,
    \end{align}
    i.e., \ $Y_\infty$ \ has Gamma distribution with parameters
    \ $2a / \sigma_1^2$ \ and \ $2b / \sigma_1^2$, \ hence
    \[
      \EE(Y_\infty^\kappa)
      = \frac{\Gamma\bigl(\frac{2a}{\sigma_1^2} + \kappa\bigr)}
             {\bigl(\frac{2b}{\sigma_1^2}\bigr)^\kappa \,
              \Gamma\bigl(\frac{2a}{\sigma_1^2}\bigr)} , \qquad
      \kappa \in \biggl( -\frac{2a}{\sigma_1^2}, \infty \biggr) .
    \]
   Especially, \ $\EE(Y_\infty) = \frac{a}{b}$.
   \ Further, if \ $a \in \bigl( \frac{\sigma_1^2}{2}, \infty \bigr)$,
    \ then \ $\EE\bigl(\frac{1}{Y_\infty}\bigr) = \frac{2b}{2a - \sigma_1^2}$.
 \item
  For all Borel measurable functions \ $f : \RR \to \RR$ \ such that
   \ $\EE(|f(Y_\infty)|) < \infty$, \ we have
   \begin{equation}\label{ergodic}
    \frac{1}{T} \int_0^T f(Y_s) \, \dd s \as \EE(f(Y_\infty)) \qquad
    \text{as \ $T \to \infty$.}
   \end{equation}
\end{enumerate}
\end{Thm}

\section{Radon--Nikodym derivatives for Heston models}
\label{section_RN}

From this section, we will consider the Heston model \eqref{Heston_SDE} with
 fixed \ $\sigma_1, \sigma_2 \in \RR_{++}$, \ $\varrho \in (-1, 1)$, \ and fixed
 initial value \ $(Y_0, X_0) = (y_0, x_0) \in \RR_{++} \times \RR$, \ and we
 will consider
 \ $\btheta := (a, \alpha, b, \beta) \in \RR_{++} \times \RR^3 =: \Theta$ \ as
 a parameter.
Note that \ $\Theta \subset \RR^4$ \ is an open subset.

Let \ $\PP_\btheta$ \ denote the probability measure induced by
 \ $(Y_t, X_t)_{t\in\RR_+}$ \ on the measurable space
 \ $(C(\RR_+, \RR^2), \cB(C(\RR_+, \RR^2)))$ \ endowed with the natural
 filtration \ $(\cG_t)_{t\in\RR_+}$, \ given by
 \ $\cG_t := \varphi_t^{-1}(\cB(C(\RR_+, \RR^2)))$, \ $t \in \RR_+$, \ where
 \ $\varphi_t : C(\RR_+, \RR^2) \to C(\RR_+, \RR^2)$ \ is the mapping
 \ $\varphi_t(f)(s) := f(t \land s)$, \ $s, t \in \RR_+$,
 \ $f\in C(\RR_+, \RR^2)$.   
\ Here \ $C(\RR_+, \RR^2)$ \ denotes the set of \ $\RR^2$-valued continuous
 functions defined on \ $\RR_+$, \ and \ $\cB(C(\RR_+, \RR^2))$ \ is the Borel
 \ $\sigma$-algebra on it.
Further, for all \ $T \in \RR_{++}$, \ let
 \ $\PP_{\btheta,T} := \PP_\btheta|_{\cG_T}$ \ be the restriction of \ $\PP_\btheta$
 \ to \ $\cG_T$.

Let us write the Heston model \eqref{Heston_SDE} in the form
 \begin{equation}\label{Heston_SDE_matrix}
  \begin{bmatrix} \dd Y_t \\ \dd X_t \end{bmatrix}
  = \begin{bmatrix} a - b Y_t \\ \alpha - \beta Y_t \end{bmatrix} \dd t
    + \sqrt{Y_t}
      \begin{bmatrix}
       \sigma_1 & 0 \\
       \sigma_2 \varrho & \sigma_2 \sqrt{1 - \varrho^2}
      \end{bmatrix}
      \begin{bmatrix} \dd W_t \\ \dd B_t \end{bmatrix} .
 \end{equation}
In order to calculate Radon--Nikodym derivatives
 \ $\frac{\dd \PP_{\bttheta,T}}{\dd \PP_{\btheta,T}}$ \ for certain
 \ $\btheta, \bttheta \in \Theta$, \ we need the following statement, which
 can be derived from formula (7.139) in Section 7.6.4 of Liptser and Shiryaev
 \cite{LipShiI}, see Barczy and Pap \cite[Lemma 3.1]{BarPap}.

\begin{Lem}\label{RN}
Let \ $a, \ta \in \bigl[ \frac{\sigma_1^2}{2}, \infty \bigr)$ \ and
 \ $b, \tb, \alpha, \talpha, \beta, \tbeta \in \RR$.
\ Let \ $\btheta := (a, \alpha, b, \beta)$ \ and
 \ $\bttheta := (\ta, \talpha, \tb, \tbeta)$.
\ Then for all \ $T \in \RR_{++}$, \ the measures \ $\PP_{\btheta,T}$ \ and
 \ $\PP_{\bttheta,T}$ \ are absolutely continuous with respect to each other,
 and
 \begin{align*}
  \log \frac{\dd \PP_{\bttheta,T}}{\dd \PP_{\btheta,T}}(Y, X)
  &=\int_0^T
     \frac{1}{Y_s}
     \begin{bmatrix}
      (\ta - \tb Y_s) - (a - b Y_s) \\
      (\talpha - \tbeta Y_s) - (\alpha - \beta Y_s)
     \end{bmatrix}^\top
     \bS^{-1}
     \begin{bmatrix} \dd Y_s \\ \dd X_s \end{bmatrix} \\
  &\quad
    -\frac{1}{2}
     \int_0^T
      \frac{1}{Y_s}
      \begin{bmatrix}
       (\ta - \tb Y_s) - (a - b Y_s) \\
       (\talpha - \tbeta Y_s) - (\alpha - \beta Y_s)
      \end{bmatrix}^\top
      \bS^{-1}
      \begin{bmatrix}
       (\ta - \tb Y_s) + (a - b Y_s) \\
       (\talpha - \tbeta Y_s) + (\alpha - \beta Y_s)
      \end{bmatrix}
      \dd s ,
 \end{align*}
 where
 \begin{equation}\label{bSigma}
  \bS := \begin{bmatrix}
          \sigma_1 & 0 \\
          \sigma_2 \varrho & \sigma_2 \sqrt{1 - \varrho^2}
         \end{bmatrix}
         \begin{bmatrix}
          \sigma_1 & \sigma_2 \varrho \\
          0 & \sigma_2 \sqrt{1 - \varrho^2}
         \end{bmatrix}
       = \begin{bmatrix}
          \sigma_1^2 & \varrho \sigma_1 \sigma_2 \\
          \varrho \sigma_1 \sigma_2 & \sigma_2^2
         \end{bmatrix} .
 \end{equation}
Moreover, the process
 \begin{equation}\label{RNmart}
  \left( \frac{\dd \PP_{\bttheta,T}}{\dd \PP_{\btheta,T}} \right)_{T\in\RR_+}
 \end{equation}
 is a \ $\PP_\btheta$-martingale with respect to the filtration
 \ $(\cG_T)_{T\in\RR_+}$.
\end{Lem}

The martingale property of the process \eqref{RNmart} is a consequence of
 Theorem 3.4 in Chapter III of Jacod and Shiryaev \cite{JSh}.

In order to investigate convergence of the family
 \begin{equation}\label{cET}
  (\cE_T)_{T\in\RR_{++}}
  := \big( C(\RR_+, \RR^2), \cB(C(\RR_+, \RR^2)) ,
           \{\PP_{\btheta,T}
             : \btheta \in \RR_{++} \times \RR^3\} \big)_{T\in\RR_{++}}
 \end{equation}
 of statistical experiments, we derive the following corollary.

\begin{Cor}\label{RN_Cor}
Let \ $a \in \bigl[ \frac{\sigma_1^2}{2}, \infty \bigr)$,
 \ $b, \alpha, \beta, \in \RR$ \ and \ $T \in \RR_{++}$.
\ Put \ $\btheta := (a, \alpha, b, \beta)$.
\ If
 \[
   \br_{\btheta,T}
   = \begin{bmatrix}
      r_{\btheta,T,1} & 0 & 0 & 0 \\
      0 & r_{\btheta,T,2} & 0 & 0 \\
      0 & 0 & r_{\btheta,T,3} & 0 \\
      0 & 0 & 0 & r_{\btheta,T,4}
     \end{bmatrix} \in \RR^{4\times4} , \qquad 
   \bh_T = \begin{bmatrix} h_{T,1} \\ h_{T,2} \\ h_{T,3} \\ h_{T,4} \end{bmatrix}
         \in \RR^4
 \]
 such that
 \ $a + r_{\btheta,T,1} h_{T,1} \in \left[ \frac{\sigma_1^2}{2}, \infty \right)$,
 \ then
 \[
   \log \frac{\dd \PP_{\btheta+\br_{\btheta,T}\bh_T,T}}{\dd \PP_{\btheta,T}}(Y, X)
   = \bh_T^\top \bDelta_{\btheta,T}(Y, X)
     - \frac{1}{2} \bh_T^\top  \bJ_{\btheta,T}(Y, X) \, \bh_T ,
 \]
 where
 \[
   \bDelta_{\btheta,T}(Y, X)
   := \br_{\btheta,T} 
      \left(\bI_2
            \otimes
            \begin{bmatrix}
             \sigma_1 & \sigma_2 \varrho \\
             0 & \sigma_2 \sqrt{1 - \varrho^2}
            \end{bmatrix}^{-1}\right)
      \begin{bmatrix}
       \int_0^T \frac{\dd W_s}{\sqrt{Y_s}} \\[2mm]
       \int_0^T \frac{\dd B_s}{\sqrt{Y_s}} \\[2mm]
       - \int_0^T \sqrt{Y_s} \, \dd W_s \\[2mm]
       - \int_0^T \sqrt{Y_s} \, \dd B_s
      \end{bmatrix} 
 \]
 and
 \[
   \bJ_{\btheta,T}(Y, X)
   := \br_{\btheta,T}
      \left( \begin{bmatrix}
              \int_0^T \frac{\dd s}{Y_s} & - T \\
              - T & \int_0^T Y_s \, \dd s
             \end{bmatrix}
             \otimes
             \bS^{-1} \right)
      \br_{\btheta,T} ,
 \]
 where \ $\bA \otimes \bB$ \ denotes the Kronecker product of matrices \ $\bA$
 \ and \ $\bB$.
\ Consequently, by Remark \ref{oO}, the quadratic approximation \eqref{LAQ} is
 valid.
\end{Cor}

\noindent{\bf Proof.}
Using equations \eqref{Heston_SDE_matrix}, we get
 \begin{align*}
  \log \frac{\dd \PP_{\bttheta,T}}{\dd \PP_{\btheta,T}}(Y, X)
  &=\int_0^T
     \frac{1}{\sqrt{Y_s}}
     \begin{bmatrix}
      (\ta - a) - (\tb - b) Y_s \\
      (\talpha - \alpha) - (\tbeta - \beta) Y_s
     \end{bmatrix}^\top
     \begin{bmatrix}
      \sigma_1 & \sigma_2 \varrho \\
      0 & \sigma_2 \sqrt{1 - \varrho^2}
     \end{bmatrix}^{-1}
     \begin{bmatrix}
      \dd W_s \\
      \dd B_s
     \end{bmatrix} \\
  &\quad
    -\frac{1}{2}
     \int_0^T
      \frac{1}{Y_s}
      \begin{bmatrix}
       (\ta - a) - (\tb - b) Y_s \\
       (\talpha - \alpha) - (\tbeta - \beta) Y_s
      \end{bmatrix}^\top
      \bS^{-1}
      \begin{bmatrix}
       (\ta - a) - (\tb - b) Y_s \\
       (\talpha - \alpha) - (\tbeta - \beta) Y_s
      \end{bmatrix}
      \dd s .
 \end{align*}
Writing \ $\br = \br_{\btheta,T}$ \ and \ $\bh = \bh_T$ \ for the sake of
 simplicity, we obtain
 \[
   \log \frac{\dd \PP_{\btheta+\br\bh,T}}{\dd \PP_{\btheta,T}}(Y, X)
   = I_1 - \frac{1}{2} I_2 ,
 \]
 where
 \begin{align*}
  I_1
  &:=\int_0^T
      \frac{1}{\sqrt{Y_s}}
      \begin{bmatrix}
       r_1 h_1 - r_3 h_3 Y_s \\
       r_2 h_2 - r_4 h_4 Y_s
      \end{bmatrix}^\top
      \begin{bmatrix}
       \sigma_1 & \sigma_2 \varrho \\
       0 & \sigma_2 \sqrt{1 - \varrho^2}
      \end{bmatrix}^{-1}
      \begin{bmatrix}
       \dd W_s \\
       \dd B_s
      \end{bmatrix} , \\
  I_2
  &:=\int_0^T
      \frac{1}{Y_s}
      \begin{bmatrix}
       r_1 h_1 - r_3 h_3 Y_s \\
       r_2 h_2 - r_4 h_4 Y_s
      \end{bmatrix}^\top
      \bS^{-1}
      \begin{bmatrix}
       r_1 h_1 - r_3 h_3 Y_s \\
       r_2 h_2 - r_4 h_4 Y_s
      \end{bmatrix}
      \dd s .
 \end{align*}
We have
 \[
   \frac{1}{\sqrt{Y_s}}
   \begin{bmatrix}
    r_1 h_1 - r_3 h_3 Y_s \\
    r_2 h_2 - r_4 h_4 Y_s
   \end{bmatrix}^\top
   = \frac{1}{\sqrt{Y_s}}
     \begin{bmatrix} h_1 \\ h_2 \end{bmatrix}^\top
     \begin{bmatrix} r_1 & 0 \\ 0 & r_2 \end{bmatrix}
     - \sqrt{Y_s}
       \begin{bmatrix} h_3 \\ h_4 \end{bmatrix}^\top
       \begin{bmatrix} r_3 & 0 \\ 0 & r_4 \end{bmatrix} ,
 \]
 hence
 \begin{align*}
  I_1
  &= \begin{bmatrix} h_1 \\ h_2 \end{bmatrix}^\top
     \begin{bmatrix} r_1 & 0 \\ 0 & r_2 \end{bmatrix}
     \begin{bmatrix}
      \sigma_1 & \sigma_2 \varrho \\
      0 & \sigma_2 \sqrt{1 - \varrho^2}
     \end{bmatrix}^{-1}
     \begin{bmatrix}
      \int_0^T \frac{\dd W_s}{\sqrt{Y_s}} \\[2mm]
      \int_0^T \frac{\dd B_s}{\sqrt{Y_s}}
     \end{bmatrix} \\
  &\quad    
     + \begin{bmatrix} h_3 \\ h_4 \end{bmatrix}^\top
       \begin{bmatrix} r_3 & 0 \\ 0 & r_4 \end{bmatrix}
       \begin{bmatrix}
        \sigma_1 & \sigma_2 \varrho \\
        0 & \sigma_2 \sqrt{1 - \varrho^2}
       \end{bmatrix}^{-1}
       \begin{bmatrix}
        - \int_0^T \sqrt{Y_s} \, \dd W_s \\[2mm]
        - \int_0^T \sqrt{Y_s} \, \dd B_s
       \end{bmatrix} 
   = \bDelta_{\btheta,T}(Y, X) ,
 \end{align*}
 and
 \begin{align*}
  I_2
  &= \int_0^T \frac{\dd s}{Y_s}
     \begin{bmatrix} h_1 \\ h_2 \end{bmatrix}^\top
     \begin{bmatrix} r_1 & 0 \\ 0 & r_2 \end{bmatrix}
     \bS^{-1}
     \begin{bmatrix} r_1 & 0 \\ 0 & r_2 \end{bmatrix}
     \begin{bmatrix} h_1 \\ h_2 \end{bmatrix}
     - T \begin{bmatrix} h_1 \\ h_2 \end{bmatrix}^\top
         \begin{bmatrix} r_1 & 0 \\ 0 & r_2 \end{bmatrix}
         \bS^{-1}
         \begin{bmatrix} r_3 & 0 \\ 0 & r_4 \end{bmatrix}
         \begin{bmatrix} h_3 \\ h_4 \end{bmatrix} \\
  &\quad    
     + T \begin{bmatrix} h_3 \\ h_4 \end{bmatrix}^\top
         \begin{bmatrix} r_3 & 0 \\ 0 & r_4 \end{bmatrix}
         \bS^{-1}
         \begin{bmatrix} r_1 & 0 \\ 0 & r_2 \end{bmatrix}
         \begin{bmatrix} h_1 \\ h_2 \end{bmatrix}
     + \int_0^T Y_s \, \dd s
       \begin{bmatrix} h_3 \\ h_4 \end{bmatrix}^\top
       \begin{bmatrix} r_3 & 0 \\ 0 & r_4 \end{bmatrix}
       \bS^{-1}
       \begin{bmatrix} r_3 & 0 \\ 0 & r_4 \end{bmatrix}
       \begin{bmatrix} h_3 \\ h_4 \end{bmatrix} \\
  &= \bJ_{\btheta,T}(Y, X) ,
 \end{align*}
 hence we conclude the assertion.
\proofend

\section{Subcritical case}
\label{section_subcrit}

\begin{Thm}\label{Thm_subcrit}
If \ $a \in \bigl( \frac{\sigma_1^2}{2}, \infty \bigr)$, \ $b \in \RR_{++}$,
 \ and \ $\alpha, \beta \in \RR$, \ then the family \ $(\cE_T)_{T\in\RR_{++}}$
 \ of statistical experiments, given in \eqref{cET}, is LAN at
 \ $\btheta := (a, \alpha, b, \beta)$ \ with scaling matrices
 \ $\br_{\btheta,T} := \frac{1}{\sqrt{T}} \bI_4$, \ $T \in \RR_{++}$, \ and with
 information matrix
 \[
   \bJ_\btheta
   := \begin{bmatrix}
       \EE\left(\frac{1}{Y_\infty}\right) & - 1 \\
       - 1 & \EE(Y_\infty)
      \end{bmatrix}
      \otimes
      \bS^{-1} . 
 \]
Consequently, the family
 \ $(C(\RR_+, \RR^4), \cB(C(\RR_+, \RR^4)),
     \{\PP_{\btheta+\bh/\sqrt{T},T} : \bh \in \RR^4\})_{T\in\RR_{++}}$
 \ of statistical experiments converges to the statistical experiment
 \ $(\RR^4 \times \RR^{4\times4}, \cB(\RR^4 \times \RR^{4\times4}),
     \{\cN_4(\bJ_\btheta \bh, \bJ_\btheta) : \bh \in \RR^4\})$
 \ as \ $T \to \infty$.
\end{Thm}

\noindent{\bf Proof.}
By part (i) of Theorem \ref{Ergodicity}, \ $\EE(Y_\infty) = \frac{a}{b}$ \ and
 \ $\EE\left(\frac{1}{Y_\infty}\right) = \frac{2b}{2a-\sigma_1^2}$, \ and hence,
 part (ii) of Theorem \ref{Ergodicity} implies
 \begin{align}\label{SLLN}
  \frac{1}{T} \int_0^T Y_s \, \dd s \as \EE(Y_\infty) \qquad \text{and} \qquad
  \frac{1}{T} \int_0^T \frac{\dd s}{Y_s}
  \as \EE\left(\frac{1}{Y_\infty}\right) \qquad
  \text{as \ $T \to \infty$.}
 \end{align}
Thus, using \ $\br_{\btheta,T} = \frac{1}{\sqrt{T}} (\bI_2 \otimes \bI_2)$,
 \ $T \in \RR_{++}$, \ and applying the identity
 \ $(\bA \otimes \bB) (\bC \otimes \bD) = (\bA \bC) \otimes (\bB \bD)$, \ we
 obtain
 \begin{align*}
  \bJ_{\btheta,T}(Y, X)
  &= (\bI_2 \otimes \bI_2)
     \left( \begin{bmatrix}
             \frac{1}{T} \int_0^T \frac{\dd s}{Y_s} & - 1 \\
             - 1 & \frac{1}{T} \int_0^T Y_s \, \dd s
            \end{bmatrix}
            \otimes
            \bS^{-1} \right)
     (\bI_2 \otimes \bI_2) \\
  &\as \begin{bmatrix}
        \EE\left(\frac{1}{Y_\infty}\right) & - 1 \\
        - 1 & \EE(Y_\infty)
       \end{bmatrix}
       \otimes
       \bS^{-1}
   = \bJ_\btheta
  \qquad \text{as \ $T \to \infty$.}
 \end{align*}
Moreover,
 \[
   \bM_T := \begin{bmatrix}
             \int_0^T \frac{\dd W_s}{\sqrt{Y_s}} \\[1mm]
             \int_0^T \frac{\dd B_s}{\sqrt{Y_s}} \\[1mm]
             - \int_0^T \sqrt{Y_s} \, \dd W_s \\[1mm]
             - \int_0^T \sqrt{Y_s} \, \dd B_s
            \end{bmatrix} , \qquad
   T \in \RR_+ ,
 \]
 is a 4-dimensional continuous local martingale with quadratic variation
 process
 \[
   \langle \bM \rangle_T
   = \begin{bmatrix}
      \int_0^T \frac{\dd s}{Y_s} & - T \\
      - T & \int_0^T Y_s \, \dd s
     \end{bmatrix}
     \otimes
     \bI_2 , \qquad t \in \RR_+ .
 \]
By \eqref{SLLN}, we have
 \[
   \frac{1}{T} \langle \bM \rangle_T
   \as \begin{bmatrix}
        \EE\left(\frac{1}{Y_\infty}\right) & - 1 \\
        - 1 & \EE(Y_\infty)
       \end{bmatrix}
       \otimes
       \bI_2 \qquad
   \text{as \ $T \to \infty$.}
 \]
Hence, Theorem \ref{THM_Zanten} yields
 \[
   \frac{1}{\sqrt{T}} \bM_T
   \distr
   \cN_4\left( \bzero , 
               \begin{bmatrix}
                \EE\left(\frac{1}{Y_\infty}\right) & - 1 \\
                - 1 & \EE(Y_\infty)
               \end{bmatrix}
               \otimes
               \bI_2 \right) \qquad
   \text{as \ $T \to \infty$,}
 \]
 consequently, as \ $T \to \infty$, \ we have
 \begin{align*}
  &\bDelta_{\btheta,T}(Y, X)
   = \frac{1}{\sqrt{T}}
     \left( \bI_2
            \otimes
            \begin{bmatrix}
             \sigma_1 & \sigma_2 \varrho \\
             0 & \sigma_2 \sqrt{1 - \varrho^2}
            \end{bmatrix}^{-1} \right)
     \bM_T \\
  &\distr
   \left( \bI_2
          \otimes
          \begin{bmatrix}
           \sigma_1 & \sigma_2 \varrho \\
           0 & \sigma_2 \sqrt{1 - \varrho^2}
          \end{bmatrix}^{-1} \right)
   \cN_4\left( \bzero , 
               \begin{bmatrix}
                \EE\left(\frac{1}{Y_\infty}\right) & - 1 \\
                - 1 & \EE(Y_\infty)
               \end{bmatrix}
               \otimes
               \bI_2 \right) \\
  &\distre
   \cN_4\left( \bzero , 
               \left( \bI_2
               \otimes
               \begin{bmatrix}
                \sigma_1 & \sigma_2 \varrho \\
                0 & \sigma_2 \sqrt{1 - \varrho^2}
               \end{bmatrix}^{-1} \right)
               \left( \begin{bmatrix}
                       \EE\left(\frac{1}{Y_\infty}\right) & - 1 \\
                       - 1 & \EE(Y_\infty)
                      \end{bmatrix}
                      \otimes
                      \bI_2 \right)\right. \\
  &\phantom{\distre\cN_4\Biggl( \bzero \quad}
             \left.
             \times
             \left( \bI_2
             \otimes
             \begin{bmatrix}
              \sigma_1 & \sigma_2 \varrho \\
              0 & \sigma_2 \sqrt{1 - \varrho^2}
             \end{bmatrix}^{-1} \right)^\top \right) \\
  &= \cN_4\left( \bzero , 
                 \begin{bmatrix}
                  \EE\left(\frac{1}{Y_\infty}\right) & - 1 \\
                  - 1 & \EE(Y_\infty)
                 \end{bmatrix}
                 \otimes
                 \bS^{-1} \right)
   = \cN_4(\bzero, \bJ_\btheta) .
 \end{align*}
Thus,
 \[
   \cL\bigl((\bDelta_{\btheta,T}, \bJ_{\btheta,T}) \mid \PP_{\btheta,T}\bigr)
   \weak \cN_4(\bzero, \bJ_\btheta) \times \delta_{\bJ_\btheta}
   \qquad \text{as \ $T \to \infty$,}
 \]
 yielding by Remark \ref{oO}, that the family \ $(\cE_T)_{T\in\RR_{++}}$ \ of
 statistical experiments is LAN at \ $\btheta$.
\proofend

\begin{Rem}\label{AOT_subcrit}
Applying Theorem \ref{AOT_LAN} for the functions
 \ $\psi_1(a, \alpha, b, \beta) := a - a_0$,
 \ $\psi_2(a, \alpha, b, \beta) := \alpha - \alpha_0$,
 \ $\psi_3(a, \alpha, b, \beta) := b - b_0$, \ and
 \ $\psi_4(a, \alpha, b, \beta) := \beta - \beta_0$,
 \ $(a, \alpha, b, \beta) \in \RR_{++} \times \RR^3$, \ we obtain that the
 family of tests that reject for values
 \begin{align*}
  S_{\btheta_0,T}^{(1)}
  &:= \frac{\sqrt{2 a_0 - \sigma_1^2}}{\sigma_1^2 \sqrt{a_0 b_0 T}}
      \int_0^T
       \frac{a_0 - b_0 Y_s}{Y_s} \, [\dd Y_s - (a_0 - b_0 Y_s) \, \dd s] , \\
  S_{\btheta_0,T}^{(2)}
  &:= \frac{\sqrt{2 a_0 - \sigma_1^2}}{\sigma_1 \sigma_2 \sqrt{a_0 b_0 T}}
      \int_0^T
       \frac{a_0 - b_0 Y_s}{Y_s} \,
       [\dd X_s - (\alpha_0 - \beta_0 Y_s) \, \dd s] , \\
  S_{\btheta_0,T}^{(3)}
  &:= \frac{1}{\sigma_1^2 \sqrt{2 b_0 T}}
      \int_0^T
       \frac{2 a_0 - \sigma_1^2 - 2 b_0 Y_s}{Y_s} \,
       [\dd Y_s - (a_0 - b_0 Y_s) \, \dd s] , \\
  S_{\btheta_0,T}^{(4)}
  &:= \frac{1}{\sigma_1 \sigma_2 \sqrt{2 b_0 T}}
      \int_0^T
       \frac{2 a_0 - \sigma_1^2 - 2 b_0 Y_s}{Y_s} \,
       [\dd X_s - (\alpha_0 - \beta_0 Y_s) \, \dd s] ,
 \end{align*}
 exceeding \ $z_\alpha$, \ respectively, are asymptotically optimal for testing
 \ $H_0^{(1)} : a \leq a_0$ \ against \ $H_1^{(1)} : a > a_0$,
 \ $H_0^{(2)} : \alpha \leq \alpha_0$ \ against
 \ $H_1^{(2)} : \alpha > \alpha_0$,
 \ $H_0^{(3)} : b \leq b_0$ \ against \ $H_1^{(3)} : b > b_0$, \ and
 \ $H_0^{(4)} : \beta \leq \beta_0$ \ against \ $H_1^{(4)} : \beta > \beta_0$,
 \ respectively, where \ $\btheta_0 = (a_0, \alpha_0, b_0, \beta_0)$ \ with
 \ $a_0 \in \bigl( \frac{\sigma_1^2}{2}, \infty \bigr)$, \ $b_0 \in \RR_{++}$,
 \ $\alpha_0, \beta_0 \in \RR$.
\ Indeed,
 \[
   \bJ_{\btheta_0}^{-1}
   = \begin{bmatrix}
      \EE\left(\frac{1}{Y_\infty}\right) & - 1 \\
      - 1 & \EE(Y_\infty)
     \end{bmatrix}^{-1}
     \otimes
     \bS
   = \frac{2 a_0 - \sigma_1^2}{\sigma_1^2}
     \begin{bmatrix}
      \frac{a_0}{b_0} & 1 \\
      1 & \frac{2b_0}{2 a_0 - \sigma_1^2}
     \end{bmatrix}
     \otimes
     \bS ,
 \]
 hence
 \begin{align*}
  \bJ_{\btheta_0}^{-1} \bDelta_{\btheta_0,T}
  &= \left(\begin{bmatrix}
            \EE\left(\frac{1}{Y_\infty}\right) & - 1 \\
            - 1 & \EE(Y_\infty)
           \end{bmatrix}^{-1}
           \otimes
           \bS\right)
     \frac{1}{\sqrt{T}}
     \left( \bI_2
            \otimes
            \begin{bmatrix}
             \sigma_1 & \sigma_2 \varrho \\
             0 & \sigma_2 \sqrt{1 - \varrho^2}
            \end{bmatrix}^{-1} \right)
     \bM_T \\
  &= \frac{1}{\sqrt{T}}
     \left(\begin{bmatrix}
            \EE\left(\frac{1}{Y_\infty}\right) & - 1 \\
            - 1 & \EE(Y_\infty)
           \end{bmatrix}^{-1}
           \otimes
           \begin{bmatrix}
            \sigma_1 & 0 \\
            \sigma_2 \varrho & \sigma_2 \sqrt{1 - \varrho^2}
           \end{bmatrix}\right)
     \int_0^T
      \left(\begin{bmatrix}
             \frac{1}{\sqrt{Y_s}} \\
             - \sqrt{Y_s}
            \end{bmatrix}
            \otimes
            \begin{bmatrix}
             \dd W_s \\
             \dd B_s
            \end{bmatrix}\right) \\
  &= \frac{1}{\sqrt{T}}
     \int_0^T
      \left(\begin{bmatrix}
             \EE\left(\frac{1}{Y_\infty}\right) & - 1 \\
             - 1 & \EE(Y_\infty)
            \end{bmatrix}^{-1}
            \begin{bmatrix}
             \frac{1}{\sqrt{Y_s}} \\
             - \sqrt{Y_s}
            \end{bmatrix}\right)
      \otimes
      \left(\begin{bmatrix}
             \sigma_1 & 0 \\
             \sigma_2 \varrho & \sigma_2 \sqrt{1 - \varrho^2}
            \end{bmatrix}
            \begin{bmatrix}
             \dd W_s \\
             \dd B_s
            \end{bmatrix}\right) \\
  &= \frac{1}{\sqrt{T}}
     \int_0^T
      \left(\begin{bmatrix}
             \EE\left(\frac{1}{Y_\infty}\right) & - 1 \\
             - 1 & \EE(Y_\infty)
            \end{bmatrix}^{-1}
            \begin{bmatrix}
             \frac{1}{Y_s} \\
             - 1
            \end{bmatrix}\right)
      \otimes
      \begin{bmatrix}
       \dd Y_s - (a_0 - b_0 Y_s) \, \dd s \\
       \dd X_s - (\alpha_0 - \beta_0 Y_s) \, \dd s
      \end{bmatrix} \\
  &= \frac{1}{\sqrt{T}}
     \int_0^T
      \begin{bmatrix}
       \frac{(2 a_0 - \sigma_1^2) (a_0 - b_0 Y_s)}{\sigma_1^2 b_0 Y_s} \\[2mm] 
       \frac{(2 a_0 - \sigma_1^2 - 2 b_0 Y_s)}{\sigma_1^2 Y_s}
      \end{bmatrix}
      \otimes
      \begin{bmatrix}
       \dd Y_s - (a_0 - b_0 Y_s) \, \dd s \\
       \dd X_s - (\alpha_0 - \beta_0 Y_s) \, \dd s
      \end{bmatrix}
 \end{align*}
 where we used
 \[
   \begin{bmatrix}
    \sigma_1 & 0 \\
    \sigma_2 \varrho & \sigma_2 \sqrt{1 - \varrho^2}
   \end{bmatrix}
   \begin{bmatrix}
    \dd W_s \\
    \dd B_s
   \end{bmatrix}
   = \frac{1}{\sqrt{Y_s}}
     \begin{bmatrix}
      \dd Y_s - (a_0 - b_0 Y_s) \, \dd s \\
      \dd X_s - (\alpha_0 - \beta_0 Y_s) \, \dd s
     \end{bmatrix}
 \]
 following from \eqref{Heston_SDE_matrix}, and
 \ $(\psi_i)'(a_0, \alpha_0, b_0, \beta_0) = \be_i$, \ $i \in \{1, 2, 3, 4\}$.
\end{Rem}

\section{Critical case}
\label{section_crit}

\begin{Thm}\label{Thm_crit}
If \ $a \in \bigl( \frac{\sigma_1^2}{2}, \infty \bigr)$, \ $b = 0$, \ and
 \ $\alpha, \beta \in \RR$, \ then the family \ $(\cE_T)_{T\in\RR_{++}}$ \ of
 statistical experiments, given in \eqref{cET}, is LAQ at
 \ $\btheta := (a, \alpha, b, \beta)$ \ with scaling matrices
 \[
   \br_{\btheta,T}
   := \begin{bmatrix}
       \frac{1}{\sqrt{\log T}} & 0 \\
       0 & \frac{1}{T}
      \end{bmatrix}
      \otimes
      \bI_2 , \quad T \in \RR_{++} , \qquad
 \]
 and with
 \begin{equation}\label{aabb}
  \bigl(\bDelta_{\btheta,T}(Y, X) , \bJ_{\btheta,T}(Y, X)\bigr)
  \distr (\bDelta_\btheta, \bJ_\btheta) \qquad \text{as \ $T \to \infty$,}
 \end{equation}
 where
 \[
   \bDelta_\btheta
   := \begin{bmatrix}
       \left(a-\frac{\sigma_1^2}{2}\right)^{-1/2}
       \begin{bmatrix}
        \sigma_1 & \sigma_2 \varrho \\
        0 & \sigma_2 \sqrt{1 - \varrho^2}
       \end{bmatrix}^{-1}
       \bZ_2 \\
       \bS^{-1}
       \begin{bmatrix}
        a - \cY_1 \\
        \alpha - \cX_1
       \end{bmatrix}
      \end{bmatrix} , \qquad
   \bJ_\btheta
   := \begin{bmatrix}
       \left(a-\frac{\sigma_1^2}{2}\right)^{-1} & 0 \\
       0 & \int_0^1 \cY_s \, \dd s
      \end{bmatrix}
      \otimes
      \bS^{-1} ,
 \]
 where \ $(\cY_t, \cX_t)_{t\in\RR_+}$ \ is the unique strong solution of the SDE
 \begin{align}\label{help_limit_YX}
  \begin{cases}
   \dd \cY_t = a \, \dd t + \sigma_1 \sqrt{\cY_t} \, \dd \cW_t , \\
   \dd \cX_t = \alpha \, \dd t
               + \sigma_2  \sqrt{\cY_t}
                 \bigl(\varrho \, \dd \cW_t
                       + \sqrt{1 - \varrho^2} \, \dd \cB_t\bigr) ,
  \end{cases} \qquad t \in \RR_+ ,
 \end{align}
 with initial value \ $(\cY_0, \cX_0) = (0, 0)$, \ where
 \ $(\cW_t, \cB_t)_{t\in\RR_+}$ \ is a $2$-dimensional standard Wiener process,
 \ $\bZ_2$ \ is a $2$-dimensional standard normally distributed random vector
 independent of \ $\bigl(\cY_1, \int_0^1 \cY_t \, \dd t, \cX_1\bigr)$, \ and
 \ $\bS$ \ is defined in \eqref{bSigma}.
Consequently, the family
 \ $(C(\RR_+, \RR^4), \cB(C(\RR_+, \RR^4)),
     \{\PP_{\btheta+\br_{\btheta,T}\bh,T} : \bh \in \RR^4\})_{T\in\RR_{++}}$
 \ of statistical experiments converges to the statistical experiment
 \ $(\RR^4 \times \RR^{4\times4}, \cB(\RR^4 \times \RR^{4\times4}),
     \{\QQ_{\btheta,\bh} : \bh \in \RR^4\})$
 \ as \ $T \to \infty$, \ where
 \[
   \QQ_{\btheta,\bh}(B)
   := \EE\left( \exp\left\{ \bh^\top \bDelta_\btheta
                            - \frac{1}{2} \bh^\top \bJ_\btheta \bh \right\}
                \bbone_B(\bDelta_\btheta, \bJ_\btheta) \right) , \qquad
   B \in \cB(\RR^p \times \RR^{4\times4}) , \quad \bh \in \RR^4 .
 \]
If \ $b = 0$ \ and \ $\beta \in \RR$ \ are fixed, then the subfamily
 \[
   \Bigl( C(\RR_+, \RR^2), \cB(C(\RR_+, \RR^2)) ,
          \Bigl\{\PP_{\btheta,T}
          : a \in \Bigl( \frac{\sigma_1^2}{2}, \infty \Bigr) , \,
            \alpha \in \RR \Bigr\} \Bigr)_{T\in\RR_{++}}
 \]
 of statistical experiments is LAN at \ $(a, \alpha)$ \ with scaling matrices
 \ $\br^{(1)}_{\btheta,T} := \frac{1}{\sqrt{\log T}} \bI_2$, \ $T \in \RR_{++}$,
 \ and with information matrix
 \ $\bJ^{(1)}_\btheta := \left(a-\frac{\sigma_1^2}{2}\right)^{-1} \bS^{-1}$.
\ Consequently, the family
 \ $(C(\RR_+, \RR^2), \cB(C(\RR_+, \RR^2)),
     \{\PP_{\btheta+\bh/\sqrt{\log T},T} : \bh_1 \in \RR^2\})_{T\in\RR_{++}}$
 \ of statistical experiments converges to the statistical experiment
 \ $(\RR^2 \times \RR^{2\times2}, \cB(\RR^2 \times \RR^{2\times2}),
     \{\cN_4(\bJ_\btheta^{(1)} \bh_1, \bJ_\btheta^{(1)}) : \bh_1 \in \RR^2\})$
 \ as \ $T \to \infty$, \ where \ $\bh := (\bh_1, \bzero)^\top \in \RR^4$.
\end{Thm}

\noindent{\bf Proof.}
We have
 \begin{align*}
  &\bDelta_{\btheta,T}(Y, X)
   = \left( \begin{bmatrix}
             \frac{1}{\sqrt{\log T}} & 0 \\
             0 & \frac{1}{T}
            \end{bmatrix}
            \otimes
            \bI_2 \right)
     \left( \bI_2
            \otimes
            \begin{bmatrix}
             \sigma_1 & \sigma_2 \varrho \\
             0 & \sigma_2 \sqrt{1 - \varrho^2}
            \end{bmatrix}^{-1} \right)
     \begin{bmatrix}
      \int_0^T \frac{\dd W_s}{\sqrt{Y_s}} \\[2mm]
      \int_0^T \frac{\dd B_s}{\sqrt{Y_s}} \\[2mm]
      - \int_0^T \sqrt{Y_s} \, \dd W_s \\[2mm]
      - \int_0^T \sqrt{Y_s} \, \dd B_s
     \end{bmatrix} \\
  &= \left( \begin{bmatrix}
             \left(\frac{1}{\log T}
                   \int_0^T \frac{\dd s}{Y_s}\right)^{1/2} & 0 \\
             0 & \left(\frac{1}{T^2} \int_0^T Y_s \, \dd s \right)^{1/2}
            \end{bmatrix}
            \otimes
            \begin{bmatrix}
             \sigma_1 & \sigma_2 \varrho \\
             0 & \sigma_2 \sqrt{1 - \varrho^2}
            \end{bmatrix}^{-1} \right)
     \begin{bmatrix}
      \frac{\int_0^T \frac{\dd W_s}{\sqrt{Y_s}}}
           {\left(\int_0^T \frac{\dd s}{Y_s}\right)^{1/2}} \\[4mm]
      \frac{\int_0^T \frac{\dd B_s}{\sqrt{Y_s}}}
           {\left(\int_0^T \frac{\dd s}{Y_s}\right)^{1/2}} \\[4mm]
      - \frac{\int_0^T \sqrt{Y_s} \, \dd W_s}
             {\left(\int_0^T Y_s \, \dd s \right)^{1/2}} \\[4mm]
      - \frac{\int_0^T \sqrt{Y_s} \, \dd B_s}
             {\left(\int_0^T Y_s \, \dd s \right)^{1/2}}
     \end{bmatrix}
 \end{align*}
 and
 \begin{align*}
  \bJ_{\btheta,T}(Y, X)
  &= \left( \begin{bmatrix}
             \frac{1}{\sqrt{\log T}} & 0 \\
             0 & \frac{1}{T}
            \end{bmatrix}
            \otimes
            \bI_2 \right)
     \left( \begin{bmatrix}
             \int_0^T \frac{\dd s}{Y_s} & - T \\
             - T & \int_0^T Y_s \, \dd s
            \end{bmatrix}
            \otimes
            \bS^{-1} \right)
     \left( \begin{bmatrix}
             \frac{1}{\sqrt{\log T}} & 0 \\
             0 & \frac{1}{T}
            \end{bmatrix}
            \otimes
            \bI_2 \right) \\
  &= \begin{bmatrix}
      \frac{1}{\log T}\int_0^T \frac{\dd s}{Y_s} & - \frac{1}{\sqrt{\log T}} \\
      - \frac{1}{\sqrt{\log T}} & \frac{1}{T^2} \int_0^T Y_s \, \dd s
     \end{bmatrix}
     \otimes
     \bS^{-1} .
 \end{align*}
It is known that
 \begin{align}\label{1/Y}
  \frac{1}{\log T} \int_0^T \frac{\dd s}{Y_s}
  \stoch \left(a - \frac{\sigma_1^2}{2}\right)^{-1} \qquad
  \text{as \ $T \to \infty$,}
 \end{align}
 see, e.g., Overbeck \cite[Lemma 5]{Ove} or Ben Alaya and Kebaier
 \cite[Proposition 2]{BenKeb1}.
Consequently, \eqref{aabb} will follow from
 \begin{align}\label{aabb1}
  \begin{aligned}
   &\Biggl( \frac{\int_0^T \frac{\dd W_s}{\sqrt{Y_s}}}
                 {\bigl(\int_0^T \frac{\dd s}{Y_s}\bigr)^{1/2}} ,
            \frac{\int_0^T \frac{\dd B_s}{\sqrt{Y_s}}}
                 {\bigl(\int_0^T \frac{\dd s}{Y_s}\bigr)^{1/2}} ,
            \frac{\int_0^T \sqrt{Y_s} \, \dd W_s}
                 {\bigl(\int_0^T Y_s \, \dd s\bigr)^{1/2}} ,
            \frac{\int_0^T \sqrt{Y_s} \, \dd B_s}
                 {\bigl(\int_0^T Y_s \, \dd s\bigr)^{1/2}} ,
            \frac{Y_T}{T} ,
            \frac{1}{T^2} \int_0^T Y_s \, \dd s \Biggr) \\
  &\distr
   \biggl( \bZ_2,
           \frac{\cY_1 - a}{\sigma_1 \bigl(\int_0^1 \cY_s \, \dd s\bigr)^{1/2}},
           Z_3, \cY_1, \int_0^1 \cY_s \, \dd s \biggr)
  \end{aligned}
 \end{align}
 as \ $T \to \infty$, \ where \ $Z_3$ \ is a standard normally distributed
 random variable independent of
 \ $\bigl(\bZ_2, \cY_1, \int_0^1 \cY_s \, \dd s\bigr)$.
\ Indeed,
 \[
   \Bigl(\bDelta_{\btheta,T}\bigl((Y_s, X_s)_{s\in[0,T]}\bigr),
         \bJ_{\btheta,T}\bigl((Y_s, X_s)_{s\in[0,T]}\bigr)\Bigr)
   \distr (\widetilde{\bDelta}_\btheta, \bJ_\btheta) \qquad
   \text{as \ $T \to \infty$,}
 \]
 where, by \eqref{bSigma},
 \begin{align*}
  \widetilde{\bDelta}_\btheta
  &:= \left( \begin{bmatrix}
              \left(a-\frac{\sigma_1^2}{2}\right)^{-1/2} & 0 \\
              0 & \left(\int_0^1 \cY_s \, \dd s\right)^{1/2}
             \end{bmatrix}
             \otimes 
             \begin{bmatrix}
              \sigma_1 & \sigma_2 \varrho \\ 
              0 & \sigma_2 \sqrt{1 - \varrho^2}
             \end{bmatrix}^{-1} \right)
      \begin{bmatrix}
       \bZ_2 \\
       \frac{a - \cY_1}{\sigma_1 \left(\int_0^1 \cY_s \, \dd s\right)^{1/2}} \\
       - Z_3
      \end{bmatrix} \\
  &=  \begin{bmatrix}
       \left(a-\frac{\sigma_1^2}{2}\right)^{-1/2}
       \begin{bmatrix}
        \sigma_1 & \sigma_2 \varrho \\ 
        0 & \sigma_2 \sqrt{1 - \varrho^2}
       \end{bmatrix}^{-1}
       \bZ_2 \\
       \bS^{-1}
       \begin{bmatrix}
        a - \cY_1 \\
        \frac{\sigma_2 \varrho}{\sigma_1} (a - \cY_1) 
        - \sigma_2 \sqrt{1 - \varrho^2} Z_3
          \left(\int_0^1 \cY_s \, \dd s\right)^{1/2}
       \end{bmatrix}
      \end{bmatrix} ,
 \end{align*}
 and
 \ $(\widetilde{\bDelta}_\btheta, \bJ_\btheta)
    \distre (\bDelta_\btheta, \bJ_\btheta)$,
 \ since 
 \[
   \Biggl( \cY_1, \int_0^1 \cY_s \, \dd s,
           \frac{\sigma_2 \varrho}{\sigma_1}
           \frac{\cY_1 - a}{\int_0^1 \cY_s \, \dd s}
           + \frac{\sigma_2 \sqrt{1 - \varrho^2}}
                  {\bigl(\int_0^1 \cY_s \, \dd s\bigr)^{1/2}}
             Z_3 \Biggr)
   \distre
   \biggl( \cY_1, \int_0^1 \cY_s \, \dd s,
           \frac{\cX_1 - \alpha}{\int_0^1 \cY_s \, \dd s} \biggr) ,
 \]
 and \ $\bZ_2$ \ is independent of \ $(Z_3, \cY_1, \int_0^1 \cY_s \, \dd s)$
 \ and of \ $(\cY_1, \int_0^1 \cY_s \, \dd s, \cX_1)$, \ see Barczy and Pap
 \cite[Equation (6.9)]{BarPap}.

We will prove \eqref{aabb1} using continuity theorem.
We have
 \begin{align}\label{help11}
  \sigma_1 \int_0^T \frac{\dd W_s}{\sqrt{Y_s}}
  = \log Y_T - \log y_0
    + \biggl( \frac{\sigma_1^2}{2} - a\biggr)
      \int_0^T \frac{\dd s}{Y_s} , \qquad T \in \RR_+ ,
 \end{align}
 see Barczy and Pap \cite[Formula (6.16)]{BarPap}.
By \eqref{Heston_SDE} and by the assumption \ $b = 0$, \ we obtain
 \[
   \sigma_1 \int_0^T \sqrt{Y_s} \, \dd W_s
   = Y_T - y_0 - a T , \qquad T \in \RR_+ .
 \]
Consequently, \ $\int_0^T \frac{\dd W_s}{\sqrt{Y_s}}$ \ and
 \ $\int_0^T \sqrt{Y_s} \, \dd W_s$ \ are measurable with respect to the
 $\sigma$-algebra \ $\sigma(Y_s , s \in [0, T])$.
\ For all \ $(u_1, u_2, u_3, u_4, v_1, v_2) \in \RR^6$ \ and \ $T \in \RR_{++}$,
 \ we have
 \begin{align*}
  &\EE\Biggl( \exp\Biggl\{ \ii u_1
                           \frac{\int_0^T \frac{\dd W_s}{\sqrt{Y_s}}}
                                {\bigl(\int_0^T \frac{\dd s}{Y_s}\bigr)^{1/2}}
                           + \ii u_2
                             \frac{\int_0^T \frac{\dd B_s}{\sqrt{Y_s}}}
                                  {\bigl(\int_0^T
                                         \frac{\dd s}{Y_s}\bigr)^{1/2}}
                           + \ii u_3
                             \frac{\int_0^T \sqrt{Y_s} \, \dd W_s}
                                  {\bigl(\int_0^T Y_s \, \dd s\bigr)^{1/2}}
                           + \ii u_4
                             \frac{\int_0^T \sqrt{Y_s} \, \dd B_s}
                                  {\bigl(\int_0^T Y_s \, \dd s\bigr)^{1/2}} \\
  &\phantom{\EE\Biggl( \exp\Biggl\{}
                           + \ii v_1 \frac{1}{T} Y_T
                           + \ii v_2 \frac{1}{T^2}
                             \int_0^T Y_s \, \dd s \biggr\}
              \Bigg| \, Y_s , s \in [0, T] \Biggr) \\
  &= \exp\Biggl\{ \ii u_1
                  \frac{\int_0^T \frac{\dd W_s}{\sqrt{Y_s}}}
                  {\bigl(\int_0^T \frac{\dd s}{Y_s}\bigr)^{1/2}}
                  + \ii u_3
                    \frac{\int_0^T \sqrt{Y_s} \, \dd W_s}
                         {\bigl(\int_0^T Y_s \, \dd s\bigr)^{1/2}}
                  + \ii v_1 \frac{1}{T} Y_T
                  + \ii v_2 \frac{1}{T^2} \int_0^T Y_s \, \dd s \Biggr\} \\
  &\quad
     \times
     \EE\Biggl( \exp\Biggl\{ \ii
                             \int_0^T
                              \Biggl( \frac{u_2}
                                 {\bigl(\int_0^T \frac{\dd t}{Y_t}\bigr)^{1/2}}
                                 \frac{1}{\sqrt{Y_s}}
                                 + \frac{u_4}
                                      {\bigl(\int_0^T Y_t \, \dd t\bigr)^{1/2}}
                                   \sqrt{Y_s} \Biggr) \dd B_s \Biggr\}
              \Bigg| \, Y_s , s \in [0, T] \Biggr) \\
  &= \exp\Biggl\{ \ii u_1
                  \frac{\int_0^T \frac{\dd W_s}{\sqrt{Y_s}}}
                  {\bigl(\int_0^T \frac{\dd s}{Y_s}\bigr)^{1/2}}
                  + \ii u_3
                    \frac{\int_0^T \sqrt{Y_s} \, \dd W_s}
                         {\bigl(\int_0^T Y_s \, \dd s\bigr)^{1/2}}
                  + \ii v_1 \frac{1}{T} Y_T
                  + \ii v_2 \frac{1}{T^2} \int_0^T Y_s \, \dd s \Biggr\} \\
  &\quad
     \times
     \exp\Biggl\{ - \frac{1}{2}
                    \int_0^T
                     \Biggl(\frac{u_2^2}
                                 {\int_0^T \frac{\dd t}{Y_t}}
                            \frac{1}{Y_s}
                            + \frac{u_4^2}
                                   {\int_0^T Y_t \, \dd t}
                              Y_s
                            + \frac{2 u_2 u_4}
                                   {\bigl(\int_0^T \frac{\dd t}{Y_t}
                                          \int_0^T Y_t \, \dd t\bigr)^{1/2}}
                     \Biggr) \dd s \Biggr\} \\
  &= \exp\Biggl\{ \ii u_1
                  \frac{\int_0^T \frac{\dd W_s}{\sqrt{Y_s}}}
                  {\bigl(\int_0^T \frac{\dd s}{Y_s}\bigr)^{1/2}}
                  + \ii u_3
                    \frac{\int_0^T \sqrt{Y_s} \, \dd W_s}
                         {\bigl(\int_0^T Y_s \, \dd s\bigr)^{1/2}}
                  + \ii v_1 \frac{1}{T} Y_T
                  + \ii v_2 \frac{1}{T^2} \int_0^T Y_s \, \dd s \Biggr\} \\
  &\quad
     \times
     \exp\Biggl\{ - \frac{1}{2} (u_2^2 + u_4^2)
                  - \frac{T u_2 u_4}
                         {\bigl(\int_0^T \frac{\dd t}{Y_t}
                                \int_0^T Y_t \, \dd t\bigr)^{1/2}} \Biggr\} ,
 \end{align*}
 where we used the independence of \ $Y$ \ and \ $B$.
\ Consequently, the joint characteristic function of the random vector on the
 left hand side of \eqref{aabb1} takes the form
 \begin{align*}
  &\EE\Biggl( \exp\Biggl\{ \ii u_1
                           \frac{\int_0^T \frac{\dd W_s}{\sqrt{Y_s}}}
                                {\bigl(\int_0^T \frac{\dd s}{Y_s}\bigr)^{1/2}}
                           + \ii u_2
                             \frac{\int_0^T \frac{\dd B_s}{\sqrt{Y_s}}}
                                  {\bigl(\int_0^T \frac{\dd s}{Y_s}\bigr)^{1/2}}
                           + \ii u_3
                             \frac{\int_0^T \sqrt{Y_s} \, \dd W_s}
                                  {\bigl(\int_0^T Y_s \, \dd s\bigr)^{1/2}}
                           + \ii u_4
                             \frac{\int_0^T \sqrt{Y_s} \, \dd B_s}
                                  {\bigl(\int_0^T Y_s \, \dd s\bigr)^{1/2}} \\
  &\phantom{\EE\Biggl( \exp\Biggl\{}
                         + \ii v_1 \frac{1}{T} Y_T
                         + \ii v_2 \frac{1}{T^2} \int_0^T Y_s \, \dd s \Biggr\}
      \Biggr) \\
  &= \ee^{-(u_2^2+ u_4^2)/2}
     \EE\Biggl(\exp\Biggl\{ \xi_T(u_1, u_3, v_1, v_2)
                           - \frac{T u_2 u_4}
                                  {\bigl(\int_0^T \frac{\dd t}{Y_t}
                                         \int_0^T Y_t \, \dd t\bigr)^{1/2}}
                   \Biggr\}\Biggr) ,
 \end{align*}
 where
 \[
   \xi_T(u_1, u_3, v_1, v_2)
   := \ii u_1 \frac{\int_0^T \frac{\dd W_s}{\sqrt{Y_s}}}
                   {\bigl(\int_0^T \frac{\dd s}{Y_s}\bigr)^{1/2}}
      + \ii u_3 \frac{\int_0^T \sqrt{Y_s} \, \dd W_s}
                     {\bigl(\int_0^T Y_s \, \dd s\bigr)^{1/2}}
      + \ii v_1 \frac{1}{T} Y_T
      + \ii v_2 \frac{1}{T^2} \int_0^T Y_s \, \dd s .
 \]
Ben Alaya and Kebaier \cite[proof of Theorem 6]{BenKeb2} proved
 \[
   \Biggl( \frac{\log Y_T - \log y_0
                 + \bigl( \frac{\sigma_1^2}{2} - a \bigr)
                   \int_0^T \frac{\dd s}{Y_s}}
                {\sqrt{\log T}} ,
           \frac{Y_T}{T}, \frac{1}{T^2} \int_0^T Y_s \, \dd s \Biggr)
   \distr \left( \frac{\sigma_1}{\sqrt{a-\frac{\sigma_1^2}{2}}} Z_1 ,
                 \cY_1 , \int_0^1 \cY_s \, \dd s \right)
 \]
 as \ $T \to \infty$, \ where \ $Z_1$ \ is a $1$-dimensional standard normally
 distributed random variable independent of
 \ $\left(\cY_1, \int_0^1 \cY_t \, \dd t\right)$.
\ Using \eqref{help11} we have
 \[
   \frac{\int_0^T \frac{\dd W_s}{\sqrt{Y_s}}}
        {\bigl(\int_0^T \frac{\dd s}{Y_s}\bigr)^{1/2}}
   = \frac{\frac{1}{\sqrt{\log T}}\frac{1}{\sigma_1}
     \left( \log Y_T - \log y_0
            + \left(\frac{\sigma_1^2}{2} - a\right)
              \int_0^T\frac{\dd s}{Y_s} \right)}
          {\left(\frac{1}{\log T} \int_0^T\frac{\dd s}{Y_s}\right)^{1/2}} ,
  \qquad T \in \RR_{++} ,
 \]
 and, by \eqref{1/Y}, we conclude
 \begin{equation}\label{BK}
  \begin{aligned}
   &\Biggl( \frac{\int_0^T \frac{\dd W_s}{\sqrt{Y_s}}}
                 {\bigl(\int_0^T \frac{\dd s}{Y_s}\bigr)^{1/2}} ,
            \frac{\int_0^T \sqrt{Y_s} \, \dd W_s}
                 {\bigl(\int_0^T Y_s \, \dd s\bigr)^{1/2}} ,
            \frac{Y_T}{T}, \frac{1}{T^2} \int_0^T Y_s \, \dd s  \Biggr) \\
   &\distr \biggl( Z_1 ,
                   \frac{\cY_1 - a}
                        {\sigma_1 \left(\int_0^1 \cY_s \, \dd s\right)^{1/2}} ,
                   \cY_1 , \int_0^1 \cY_s \, \dd s \biggr) \qquad
   \text{as \ $T \to \infty$,}
  \end{aligned}
 \end{equation}
 thus we derived joint convergence of four coordinates of the left hand side
 of \eqref{aabb1}.
Hence
 \begin{equation}\label{critxiconv}
  \begin{aligned}
   &\EE(\exp\{\xi_T(u_1, u_3, v_1, v_2)\}) \\
   &\to
    \EE\biggl(\exp\biggl\{ \ii u_1 Z_1
                           + \ii u_3
                             \frac{\cY_1 - a}
                                  {\sigma_1
                                   \left(\int_0^1 \cY_s \, \dd s\right)^{1/2}}
                           + \ii v_1 \cY_1
                           + \ii v_2 \int_0^1 \cY_s \, \dd s
                  \biggr\}\biggr)
  \end{aligned}
 \end{equation}
 as \ $T \to \infty$ \ for all \ $(u_1, u_3, v_1, v_2) \in \RR^4$.
\ Using \ $|\exp\{\xi_T(u_1, u_3, v_1, v_2)\}| = 1$, \ we have
 \begin{align*}
  &\Biggl|\EE\Biggl(\exp\Biggl\{ \xi_T(u_1, u_3, v_1, v_2)
                                 - \frac{T u_2 u_4}
                                        {\bigl(\int_0^T \frac{\dd t}{Y_t}
                                              \int_0^T Y_t \, \dd t\bigr)^{1/2}}
                        \Biggr\}\Biggr)
          -\EE(\exp\{\xi_T(u_1, u_3, v_1, v_2)\})\Biggr| \\
  &\qquad
   \leq \EE\Biggl(|\exp\{\xi_T(u_1, u_3, v_1, v_2)\}|
                  \Biggl|\exp\Biggl\{-\frac{T u_2 u_4}
                                           {\bigl(\int_0^T \frac{\dd t}{Y_t}
                                              \int_0^T Y_t \, \dd t\bigr)^{1/2}}
                             \Biggr\}
                         -1\Biggr|\Biggr) \\
  &\qquad
   = \EE\Biggl(\Biggl|\exp\Biggl\{-\frac{T u_2 u_4}
                                        {\bigl(\int_0^T \frac{\dd t}{Y_t}
                                              \int_0^T Y_t \, \dd t\bigr)^{1/2}}
                          \Biggr\}
                      -1\Biggr|\Biggr)
   \to 0 \qquad \text{as \ $T \to \infty$,}
 \end{align*}
 by the moment convergence theorem (see, e.g., Stroock \cite[Lemma 2.2.1]{Str}).
Indeed, by \eqref{1/Y}, \eqref{BK}, continuous mapping theorem and Slutsky's
 lemma,
 \[
   \left| \exp\Biggl\{-\frac{T u_2 u_4}
                            {\bigl(\int_0^T \frac{\dd t}{Y_t}
                                   \int_0^T Y_t \, \dd t\bigr)^{1/2}}\Biggr\}
          - 1 \right|
   = \left| \exp\Biggl\{-\frac{u_2 u_4}
                              {\sqrt{\log T}
                               \bigl(\frac{1}{\log T}
                                     \int_0^T \frac{\dd t}{Y_t}
                                     \cdot \frac{1}{T^2}
                                     \int_0^T Y_t \, \dd t\bigr)^{1/2}}\Biggr\}
            - 1 \right|
   \stoch 0
 \]
 as \ $T \to \infty$, \ and the family
 \[
   \left\{ \left| \exp\Biggl\{-\frac{T u_2 u_4}
                                    {\bigl(\int_0^T \frac{\dd t}{Y_t}
                                           \int_0^T
                                            Y_t \, \dd t\bigr)^{1/2}}\Biggr\}
                  - 1 \right| , \, \, T \in \RR_{++} \right\}
 \]
 is uniformly integrable, since, by Cauchy--Schwarz inequality,
 \[
   \Biggl| \exp\Biggl\{-\frac{T u_2 u_4}
                             {\bigl(\int_0^T \frac{\dd t}{Y_t}
                                    \int_0^T Y_t \, \dd t\bigr)^{1/2}}\Biggr\}
          - 1\Biggr|
   \leq \exp\Biggl\{\frac{T |u_2 u_4|}
                         {\bigl(\int_0^T \frac{\dd t}{Y_t}
                                 \int_0^T Y_t \, \dd t\bigr)^{1/2}}\Biggr\}
        + 1
   \leq \exp\{|u_2 u_4|\} + 1
 \]
 for all \ $T \in \RR_{++}$.
\ Using \eqref{critxiconv}, we conclude
 \begin{align*}
  &\EE\Biggl(\exp\Biggl\{\ii u_1
                         \frac{\int_0^T \frac{\dd W_s}{\sqrt{Y_s}}}
                              {\bigl(\int_0^T \frac{\dd s}{Y_s}\bigr)^{1/2}}
                         + \ii u_2
                           \frac{\int_0^T \frac{\dd B_s}{\sqrt{Y_s}}}
                                {\bigl(\int_0^T \frac{\dd s}{Y_s}\bigr)^{1/2}}
                         + \ii u_3
                           \frac{\int_0^T \sqrt{Y_s} \, \dd W_s}
                                {\bigl(\int_0^T Y_s \, \dd s\bigr)^{1/2}}
                         + \ii u_4
                           \frac{\int_0^T \sqrt{Y_s} \, \dd B_s}
                                {\bigl(\int_0^T Y_s \, \dd s\bigr)^{1/2}} \\
  &\phantom{\EE\Biggl( \exp\Biggl\{}
                         + \ii v_1 \frac{1}{T} Y_T
                         + \ii v_2 \frac{1}{T^2}
                                   \int_0^T Y_s \, \dd s \Biggr\}\Biggr) \\
  &\to \ee^{-(u_2^2+u_4^2)/2}
       \EE\biggl(\exp\Biggl\{ \ii u_1 Z_1
                              + \ii u_3
                                \frac{\cY_1 - a}
                                     {\sigma_1
                                      \left(\int_0^1 \cY_s \, \dd s\right)^{1/2}}
                              + \ii v_1 \cY_1
                              + \ii v_2 \int_0^1 \cY_s \, \dd s \Biggr\}
          \biggr)
 \end{align*}
 as \ $T \to \infty$.
\ Note that, since \ $Z_1$ \ is independent of
 \ $\bigl( \cY_1, \int_0^1 \cY_s \, \dd s \bigr)$, \ we have
 \begin{align*}
  &\ee^{-(u_2^2+u_4^2)/2}
   \EE\biggl(\exp\biggl\{ \ii u_1 Z_1
                          + \ii u_3
                            \frac{\cY_1 - a}
                                 {\sigma_1
                                  \left(\int_0^1 \cY_s \, \dd s\right)^{1/2}}
                          + \ii v_1 \cY_1
                          + \ii v_2 \int_0^1 \cY_s \, \dd s
                 \biggr\} \biggr) \\
  &= \EE(\ee^{\ii u_1 Z_1})\EE(\ee^{\ii u_2 Z_2}) \EE(\ee^{\ii u_3 Z_3})
     \EE\biggl(\exp\biggl\{ \ii u_3
                            \frac{\cY_1 - a}
                                 {\sigma_1
                                  \left(\int_0^1 \cY_s \, \dd s\right)^{1/2}}
                            + \ii v_1 \cY_1
                            + \ii v_2 \int_0^1 \cY_s \, \dd s
                   \biggr\}\biggr) ,
 \end{align*}
 where \ $(Z_2, Z_3)$ \ is a
 2-dimensional standard normally distributed random vector, independent of
 \ $\bigl( Z_1, \cY_1, \int_0^1 \cY_s \, \dd s \bigr)$, \ thus we obtain
 \eqref{aabb1} with \ $\bZ_2 := (Z_1, Z_2)$, \ and hence \eqref{aabb}, which
 yields \eqref{LAQO}.

It is known that \ $\PP(\int_0^1 \cY_s \, \dd s \in \RR_{++}) = 1$
 \ (which has been shown in the proof of Theorem 3.1 in Barczy et
 al.~\cite{BarDorLiPap}), hence \eqref{LAQJ} holds.
Finally, \eqref{LAQDJ} will follow from
 \begin{equation}\label{LAQDJcrit}
  \EE\left( \exp\left\{ \bh^\top \bDelta_\btheta
                        - \frac{1}{2} \bh^\top \bJ_\btheta \bh \right\} \right)
  = 1
 \end{equation}
 for all \ $\bh \in \RR^4$.
\ Writing \ $\bh = (\bh_1, \bh_2)^\top$, \ $\bh_1, \bh_2 \in \RR^2$, \ and using
 the independence of \ $\bZ_2$ \ and
 \ $\bigl(\cY_1, \int_0^1 \cY_t \, \dd t, \cX_1\bigr)$, \ we have
 \[
   \EE\left( \exp\left\{ \bh^\top \bDelta_\btheta
                         - \frac{1}{2} \bh^\top \bJ_\btheta \bh \right\} \right)
   = E_1 E_2 ,
 \]
 where
 \begin{align*}
  E_1 &:= \EE\left(\exp\left\{\left(a - \frac{\sigma_1^2}{2}\right)^{-1/2}
                              \bh_1^\top
                              \begin{bmatrix}
                               \sigma_1 & \sigma_2 \varrho \\
                               0 & \sigma_2 \sqrt{1 - \varrho^2}
                              \end{bmatrix}^{-1} \!\!
                              \bZ_2
                              - \frac{1}{2}
                                \left(a - \frac{\sigma_1^2}{2}\right)^{-1}
                                \bh_1^\top \bS^{-1} \bh_1\right\}\right) , \\
  E_2 &:= \EE\left(\exp\left\{\bh_2^\top \bS^{-1}
                              \begin{bmatrix}
                               a - \cY_1 \\
                               \alpha - \cX_1
                              \end{bmatrix}
                              - \frac{1}{2}
                                \left(\int_0^1 \cY_s \, \dd s\right)
                                \bh_2^\top \bS^{-1} \bh_2\right\}\right) .
 \end{align*}
The moment generating function of the 2-dimensional standard normally
 distributed random vector \ $\bZ_2$ \ has the form
 \begin{equation}\label{momgennorm}
   \EE(\ee^{\bv^\top\!\bZ_2})
   = \ee^{\|\bv\|^2/2} , \qquad \bv \in \RR^2 ,
 \end{equation}
 since
 \[
   \EE(\ee^{\bv^\top\!\bZ_2})
   = \frac{1}{2\pi} \int_{\RR^2} \ee^{\bv^\top\!\bZ_2-\|\bx\|^2/2} \, \dd \bx
   = \frac{1}{2\pi} \int_{\RR^2} \ee^{-\|\bx-\bv\|^2/2+\|\bv\|^2/2} \, \dd \bx
   = \ee^{\|\bv\|^2/2} .
 \]
Applying this with
 \[
   \bv^\top = \left(a - \frac{\sigma_1^2}{2}\right)^{-1/2} \bh_1^\top
             \begin{bmatrix}
              \sigma_1 & \sigma_2 \varrho \\
              0 & \sigma_2 \sqrt{1 - \varrho^2}
             \end{bmatrix}^{-1} , \qquad
   \|\bv\|^2 = \bv^\top \bv = \left(a - \frac{\sigma_1^2}{2}\right)^{-1}
                             \bh_1^\top \bS^{-1} \bh_1 ,
 \]
 we obtain \ $E_1 = 1$. 
\ Using Corollary \ref{RN_Cor} for the process \ $(\cY_t, \cX_t)_{t\in\RR_+}$
 \ with
 \[
   \br_{\btheta,T} = \br
   := \begin{bmatrix} 0 & 0 \\ 0 & 1 \end{bmatrix} \otimes \bI_2 , \qquad
   \bh_T = \bh
 \]
 we obtain
 \[
   \log \frac{\dd \PP_{\btheta+\br\bh,T}}{\dd \PP_{\btheta,T}}(\cY, \cX)
   = \bh^\top \bDelta_{\btheta,T}(\cY, \cX)
     - \frac{1}{2} \bh^\top  \bJ_{\btheta,T}(\cY, \cX) \, \bh , 
 \]
 where
 \begin{align*}
  &\bh^\top \bDelta_{\btheta,T}(\cY, \cX)
   = \bh^\top 
     \left( \begin{bmatrix}
             0 & 0 \\
             0 & 1
            \end{bmatrix}
            \otimes
            \bI_2 \right)
     \left( \bI_2
            \otimes
            \begin{bmatrix}
             \sigma_1 & \sigma_2 \varrho \\
             0 & \sigma_2 \sqrt{1 - \varrho^2}
            \end{bmatrix}^{-1} \right)
     \begin{bmatrix}
      \int_0^T \frac{\dd \cW_s}{\sqrt{\cY_s}} \\[2mm]
      \int_0^T \frac{\dd \cB_s}{\sqrt{\cY_s}} \\[2mm]
      - \int_0^T \sqrt{\cY_s} \, \dd \cW_s \\[2mm]
      - \int_0^T \sqrt{\cY_s} \, \dd \cB_s
     \end{bmatrix} \\
  &= - \bh_2^\top
       \begin{bmatrix}
        \sigma_1 & \sigma_2 \varrho \\
        0 & \sigma_2 \sqrt{1 - \varrho^2}
       \end{bmatrix}^{-1}
       \begin{bmatrix}
        \int_0^T \sqrt{\cY_s} \, \dd \cW_s \\[2mm]
        \int_0^T \sqrt{\cY_s} \, \dd \cB_s
       \end{bmatrix}
   = - \bh_2^\top \bS^{-1}
       \begin{bmatrix}
        \sigma_1 & 0 \\
        \sigma_2 \varrho & \sigma_2 \sqrt{1 - \varrho^2}
       \end{bmatrix}
       \begin{bmatrix}
        \int_0^T \sqrt{\cY_s} \, \dd \cW_s \\[2mm]
        \int_0^T \sqrt{\cY_s} \, \dd \cB_s
       \end{bmatrix} \\
  &= \bh_2^\top \bS^{-1}
     \begin{bmatrix}
      a T - \cY_T \\
      \alpha T - \cX_T
     \end{bmatrix}
 \end{align*}
 and
 \begin{align*}
  \bh^\top \bJ_{\btheta,T}(\cY, \cX) \bh
  &= \bh^\top 
     \left( \begin{bmatrix} 0 & 0 \\ 0 & 1 \end{bmatrix}
            \otimes
            \bI_2 \right)
     \left( \begin{bmatrix}
             \int_0^T \frac{\dd s}{\cY_s} & - T \\
             - T & \int_0^T \cY_s \, \dd s
            \end{bmatrix}
            \otimes
            \bS^{-1} \right)
     \left( \begin{bmatrix} 0 & 0 \\ 0 & 1 \end{bmatrix}
            \otimes
            \bI_2 \right)
     \bh \\
  &= \left(\int_0^T \cY_s \, \dd s\right) \bh_2^\top \bS^{-1} \bh_2 .
 \end{align*}
By Lemma \ref{RN}, the process
 \[
   \left( \frac{\dd \PP_{\btheta+\br\bh,T}}{\dd \PP_{\btheta,T}}(\cY, \cX)
   \right)_{T\in\RR_+}
   = \left( \exp\left\{ \bh_2^\top \bS^{-1}
                        \begin{bmatrix}
                         a T - \cY_T \\
                         \alpha T - \cX_T
                        \end{bmatrix}
                        - \frac{1}{2}
                          \left(\int_0^T \cY_s \, \dd s\right)
                          \bh_2^\top \bS^{-1} \bh_2 \right\} \right)_{T\in\RR_+}
 \]
 is a martingale, hence
 \[
   E_2
   = \EE\left( \frac{\dd \PP_{\btheta+\br\bh,1}}
                    {\dd \PP_{\btheta,1}}(\cY, \cX) \right)
   = \EE\left( \frac{\dd \PP_{\btheta+\br\bh,0}}
                    {\dd \PP_{\btheta,0}}(\cY, \cX) \right)
   = 1 ,
 \]
 and we conclude that the family \ $(\cE_T)_{T\in\RR{++}}$ \ of statistical
 experiments is LAQ at \ $\btheta$.
\proofend

\begin{Rem}\label{AOT_crit}
If \ $\btheta_0 = (a_0, \alpha_0, b_0, \beta_0)$ \ with
 \ $a_0 \in \bigl( \frac{\sigma_1^2}{2}, \infty \bigr)$, \ $b_0 = 0$ \ and
 \ $\alpha_0, \beta_0 \in \RR$, \ then applying Theorem \ref{AOT_LAN} for the
 functions \ $\psi_1(a, \alpha, b, \beta) := a - a_0$ \ and
 \ $\psi_2(a, \alpha, b, \beta) := \alpha - \alpha_0$,
 \ $(a, \alpha, b, \beta) \in \RR_{++} \times \RR^3$, \ we obtain that the
 family of tests that reject for values
 \begin{align*}
  S_{\btheta_0,T}^{(1)}
  &:= \frac{\sqrt{2 a_0 - \sigma_1^2}}{\sigma_1 \sqrt{2 \log T}}
      \int_0^T \frac{\dd Y_s - (a_0 - b_0 Y_s) \, \dd s}{Y_s} , \\
  S_{\btheta_0,T}^{(2)}
  &:= \frac{\sqrt{2 a_0 - \sigma_1^2}}{\sigma_2 \sqrt{2 \log T}}
      \int_0^T \frac{\dd X_s - (\alpha_0 - \beta_0 Y_s) \, \dd s}{Y_s} ,
 \end{align*}
 exceeding \ $z_\alpha$, \ respectively, are asymptotically optimal for testing
 \ $H_0^{(1)} : a \leq a_0$ \ against \ $H_1^{(1)} : a > a_0$, \ and
 \ $H_0^{(2)} : \alpha \leq \alpha_0$ \ against
 \ $H_1^{(2)} : \alpha > \alpha_0$, \ respectively.
Indeed,
 \ $\bigl(\bJ_{\btheta_0}^{(1)}\bigr)^{-1}
    = \Bigl(a_0 - \frac{\sigma_1^2}{2}\Bigr) \bS$,
 \begin{align*}
  \bDelta_{\btheta_0,T}
  &= \left(\begin{bmatrix}
            \frac{1}{\sqrt{\log T}} & 0 \\
            0 & \frac{1}{T}
           \end{bmatrix}
           \otimes
           \begin{bmatrix}
            \sigma_1 & \sigma_2 \varrho \\
            0 & \sigma_2 \sqrt{1 . \varrho^2}
           \end{bmatrix}^{-1}\right)
     \int_0^T
      \left(\begin{bmatrix}
             \frac{1}{\sqrt{Y_s}} \\
             - \sqrt{Y_s}
            \end{bmatrix}
            \otimes
            \begin{bmatrix}
             \dd W_s \\
             \dd B_s
            \end{bmatrix}\right) \\
  &= \int_0^T
      \left(\begin{bmatrix}
             \frac{1}{\sqrt{\log T}} & 0 \\
             0 & \frac{1}{T}
            \end{bmatrix}
            \begin{bmatrix}
             \frac{1}{\sqrt{Y_s}} \\
             - \sqrt{Y_s}
            \end{bmatrix}\right)
      \otimes
      \left(\begin{bmatrix}
             \sigma_1 & \sigma_2 \varrho \\
             0 & \sigma_2 \sqrt{1 . \varrho^2}
            \end{bmatrix}^{-1}
            \begin{bmatrix}
             \dd W_s \\
             \dd B_s
            \end{bmatrix}\right) \\
  &= \int_0^T
      \begin{bmatrix}
       \frac{1}{\sqrt{Y_s \log T}} \\
       - \frac{\sqrt{Y_s}}{T}
      \end{bmatrix}
      \otimes
      \left(\bS^{-1}
            \frac{1}{\sqrt{Y_s}}
            \begin{bmatrix}
             \dd Y_s - (a_0 - b_0 Y_s) \, \dd s \\
             \dd X_s - (\alpha_0 - \beta_0 Y_s) \, \dd s
            \end{bmatrix}\right) ,
 \end{align*}
 where we used
 \[
   \begin{bmatrix}
    \sigma_1 & \sigma_2 \varrho \\
    0 & \sigma_2 \sqrt{1 . \varrho^2}
   \end{bmatrix}^{-1}
   \begin{bmatrix}
    \dd W_s \\
    \dd B_s
   \end{bmatrix}
   = \bS^{-1}
     \begin{bmatrix}
      \sigma_1 & 0 \\
      \sigma_2 \varrho & \sigma_2 \sqrt{1 . \varrho^2}
     \end{bmatrix}
     \begin{bmatrix}
      \dd W_s \\
      \dd B_s
     \end{bmatrix} \\
   = \bS^{-1}
     \frac{1}{\sqrt{Y_s}}
     \begin{bmatrix}
      \dd Y_s - (a_0 - b_0 Y_s) \, \dd s \\
      \dd X_s - (\alpha_0 - \beta_0 Y_s) \, \dd s
     \end{bmatrix} ,
 \]
 following from \eqref{Heston_SDE_matrix}, thus
 \[
   \bDelta_{\btheta_0,T}^{-1}
   = \frac{1}{Y_s \sqrt{\log T}}
     \bS^{-1}
     \int_0^T
      \begin{bmatrix}
       \dd Y_s - (a_0 - b_0 Y_s) \, \dd s \\
       \dd X_s - (\alpha_0 - \beta_0 Y_s) \, \dd s
      \end{bmatrix} ,
 \]
 hence
 \[
   \bigl(\bJ_{\btheta_0}^{(1)}\bigr)^{-1} \bDelta_{\btheta_0,T}^{-1}
   = \Bigl(a_0 - \frac{\sigma_1^2}{2}\Bigr)
     \frac{1}{\sqrt{\log T}}
     \int_0^T
      \frac{1}{Y_s}
      \begin{bmatrix}
       \dd Y_s - (a_0 - b_0 Y_s) \, \dd s \\
       \dd X_s - (\alpha_0 - \beta_0 Y_s) \, \dd s
      \end{bmatrix} ,
 \]
 and \ $\psi_i'(a_0, \alpha_0, b_0, \beta_0) = \be_i$, \ $i \in \{1, 2\}$.
\end{Rem}

\section{Supercritical case}
\label{section_supercrit}

\begin{Thm}\label{Thm_supercrit}
If \ $a \in \bigl[ \frac{\sigma_1^2}{2}, \infty \bigr)$, \ $b \in \RR_{--}$,
 \ and \ $\alpha, \beta \in \RR$, \ then the family \ $(\cE_T)_{T\in\RR_{++}}$
 \ of statistical experiments, given in \eqref{cET}, is not LAQ at
 \ $\btheta := (a, \alpha, b, \beta)$ \ with scaling matrices
 \[
   \br_{\btheta,T}
   := \begin{bmatrix}
       1 & 0 \\
       0 & \ee^{bT/2}
      \end{bmatrix}
      \otimes
      \bI_2 , \quad T \in \RR_{++} , \qquad
 \]
 although
 \begin{equation}\label{aabb_super}
  \bigl(\bDelta_{\btheta,T}(Y, X) , \bJ_{\btheta,T}(Y, X)\bigr)
  \distr (\bDelta_\btheta, \bJ_\btheta) \qquad \text{as \ $T \to \infty$,}
 \end{equation}
 with
 \[
   \bDelta_\btheta
   := \left( \bI_2
             \otimes
             \begin{bmatrix}
              \sigma_1 & \sigma_2 \varrho \\  
              0 & \sigma_2 \sqrt{1 - \varrho^2}
             \end{bmatrix}^{-1} \right)
      \begin{bmatrix}
       \sigma_1^{-1} \tcV \\
       Z_1 \\
       \left(- \frac{\tcY_{-1/b}}{b}\right)^{1/2} \bZ_2
      \end{bmatrix} , \quad
   \bJ_\btheta
   := \begin{bmatrix}
       \int_0^{-1/b} \tcY_u \, \dd u & 0 \\
       0 & - \frac{\tcY_{-1/b}}{b}
      \end{bmatrix}
      \otimes
      \bS^{-1} ,
 \]
 where \ $(\tcY_t)_{t\in\RR_+}$ \ is a CIR process given by the SDE
 \[
   \dd \tcY_t = a \dd t + \sigma_1 \sqrt{\tcY_t} \, \dd \cW_t ,
   \qquad t \in \RR_+ ,
 \]
 with initial value \ $\tcY_0 = y_0$, \ where \ $(\cW_t)_{t\in\RR_+}$ \ is a
 standard Wiener process,
 \[
   \tcV
   :=\log \tcY_{-1/b} - \log y_0
     - \left(a - \frac{\sigma_1^2}{2}\right) \int_0^{-1/b} \tcY_u \, \dd u ,
 \]
 $Z_1$ \ is a $1$-dimensional standard normally distributed random variable,
 \ $\bZ_2$ \ is a $2$-dimensional standard normally distributed random vector
 such that \ $(\tcY_{-1/b},\int_0^{-1/b}\tcY_u \, \dd u)$, \ $Z_1$ \ and
 \ $\bZ_2$ \ are independent, and \ $\bS$ \ is defined in \eqref{bSigma}.
Moreover, \eqref{LAQJ} also holds, but \eqref{LAQDJ} is not valid.

If \ $a \in \bigl( \frac{\sigma_1^2}{2}, \infty \bigr)$ \ and
 \ $\alpha \in \RR$ \ are fixed, then the subfamily
 \[
   \bigl( C(\RR_+, \RR^2), \cB(C(\RR_+, \RR^2)) ,
          \bigl\{\PP_{\btheta,T}
          : b \in \RR_{--} , \,
            \beta \in \RR \bigr\} \bigr)_{T\in\RR_{++}}
 \]
 of statistical experiments is LAMN at \ $(b, \beta)$ \ with scaling matrices
 \ $\br^{(2)}_{\btheta,T} := \ee^{bT/2} \bI_2$, \ $T \in \RR_{++}$, \ and with
 \[
   \bDelta^{(2)}_\btheta := \left(- \frac{\tcY_{-1/b}}{b}\right)^{1/2}
                    \begin{bmatrix}
                     \sigma_1 & \sigma_2 \varrho \\  
                     0 & \sigma_2 \sqrt{1 - \varrho^2}
                    \end{bmatrix}^{-1}
                    \bZ_2 , \qquad
   \bJ^{(2)}_\btheta := \left(- \frac{\tcY_{-1/b}}{b}\right) \bS^{-1} .
 \]
Consequently, the family
 \ $(C(\RR_+, \RR^2), \cB(C(\RR_+, \RR^2)),
     \{\PP_{\btheta+\ee^{bT/2}\bh,T} : \bh_2 \in \RR^2\})_{T\in\RR_{++}}$
 \ of statistical experiments converges to the statistical experiment
 \ $(\RR^2 \times \RR^{2\times2}, \cB(\RR^2 \times \RR^{2\times2}),
     \{\cL((\bDelta^{(2)}_\btheta + \bJ^{(2)}_\btheta \bh_2, \bJ^{(2)}_\btheta)
           \mid \PP)
       : \bh_2 \in \RR^2\})$
 \ as \ $T \to \infty$, \ where \ $\bh := (\bzero, \bh_2)^\top \in \RR^4$.
\end{Thm}

\noindent{\bf Proof.}
We have
 \begin{align*}
  &\bDelta_{\btheta,T}(Y, X)
   = \left( \begin{bmatrix}
             1 & 0 \\
             0 & \ee^{bT/2}
            \end{bmatrix}
            \otimes
            \bI_2 \right)
     \left( \bI_2
            \otimes
            \begin{bmatrix}
             \sigma_1 & \sigma_2 \varrho \\
             0 & \sigma_2 \sqrt{1 - \varrho^2}
            \end{bmatrix}^{-1} \right)
     \begin{bmatrix}
      \int_0^T \frac{\dd W_s}{\sqrt{Y_s}} \\[2mm]
      \int_0^T \frac{\dd B_s}{\sqrt{Y_s}} \\[2mm]
      - \int_0^T \sqrt{Y_s} \, \dd W_s \\[2mm]
      - \int_0^T \sqrt{Y_s} \, \dd B_s
     \end{bmatrix} \\
  &= \left( \begin{bmatrix}
             \left(\int_0^T \frac{\dd s}{Y_s}\right)^{1/2} & 0 \\
             0 & \left(\ee^{bT} \int_0^T Y_s \, \dd s \right)^{1/2}
            \end{bmatrix}
            \otimes
            \begin{bmatrix}
             \sigma_1 & \sigma_2 \varrho \\
             0 & \sigma_2 \sqrt{1 - \varrho^2}
            \end{bmatrix}^{-1} \right)
     \begin{bmatrix}
      \frac{\int_0^T \frac{\dd W_s}{\sqrt{Y_s}}}
           {\left(\int_0^T \frac{\dd s}{Y_s}\right)^{1/2}} \\[4mm]
      \frac{\int_0^T \frac{\dd B_s}{\sqrt{Y_s}}}
           {\left(\int_0^T \frac{\dd s}{Y_s}\right)^{1/2}} \\[4mm]
      - \frac{\int_0^T \sqrt{Y_s} \, \dd W_s}
             {\left(\int_0^T Y_s \, \dd s \right)^{1/2}} \\[4mm]
      - \frac{\int_0^T \sqrt{Y_s} \, \dd B_s}
             {\left(\int_0^T Y_s \, \dd s \right)^{1/2}}
     \end{bmatrix}
 \end{align*}
 and
 \begin{align*}
  \bJ_{\btheta,T}(Y, X)
  &= \left( \begin{bmatrix}
             1 & 0 \\
             0 & \ee^{bT/2}
            \end{bmatrix}
            \otimes
            \bI_2 \right)
     \left( \begin{bmatrix}
             \int_0^T \frac{\dd s}{Y_s} & - T \\
             - T & \int_0^T Y_s \, \dd s
            \end{bmatrix}
            \otimes
            \bS^{-1} \right)
     \left( \begin{bmatrix}
             1 & 0 \\
             0 & \ee^{bT/2}
            \end{bmatrix}
            \otimes
            \bI_2 \right) \\
  &= \begin{bmatrix}
      \int_0^T \frac{\dd s}{Y_s} & - T \ee^{bT/2} \\
      - T \ee^{bT/2} & \ee^{bT} \int_0^T Y_s \, \dd s
     \end{bmatrix}
     \otimes
     \bS^{-1} .
 \end{align*}
We have
 \[
   \sigma_1 \int_0^T \frac{\dd W_s}{\sqrt{Y_s}}
   = \log Y_T - \log y_0
     + \left( \frac{\sigma_1^2}{2} - a\right)
       \int_0^T \frac{\dd s}{Y_s} + b T , \qquad T \in \RR_{++} .
 \]
 see Barczy and Pap \cite[Formula (4.10)]{BarPap}.
Moreover,
 \[
   \ee^{bT} Y_T \as V , \qquad
   \ee^{bT} \int_0^T Y_s \, \dd s \as -\frac{V}{b} , \qquad
   \int_0^T \frac{\dd s}{Y_s} \as \int_0^\infty \frac{\dd s}{Y_s} , \qquad
   \text{as \ $T \to \infty$,} 
 \]
 see Barczy and Pap \cite[Formulae (4.7) and (4.9)]{BarPap}.
Thus,
 \begin{equation}\label{a_super}
   \frac{\sigma_1 \int_0^T \frac{\dd W_s}{\sqrt{Y_s}}}
         {\int_0^T \frac{\dd s}{Y_s}}
   = \frac{\log(\ee^{bT} Y_T) - \log y_0}{\int_0^T \frac{\dd s}{Y_s}}
      + \frac{\sigma_1^2}{2} - a
   \as \frac{\log V - \log y_0}{\int_0^\infty \frac{\dd s}{Y_s}}
        + \frac{\sigma_1^2}{2} - a
 \end{equation}
 as \ $T \to \infty$.
\ By Theorem 4 in Ben Alaya and Kebaier \cite{BenKeb2},
 \[
   \left( V, \int_0^\infty \frac{\dd s}{Y_s} \right)
   \distre \left( \tcY_{-1/b} , \int_0^{-1/b} \tcY_u \, \dd u \right) ,
 \]
 hence
 \[
   \frac{\log V - \log y_0}{\int_0^\infty \frac{\dd s}{Y_s}}
   + \frac{\sigma_1^2}{2} - a
   \distre
   \frac{\log \tcY_{-1/b} - \log y_0}{\int_0^{-1/b} \tcY_u \, \dd u}
   + \frac{\sigma_1^2}{2} - a
   = \tcV .
 \]
Further,
 \[
   \Biggl( \frac{\int_0^T \frac{\dd B_s}{\sqrt{Y_s}}}
                {\bigl(\int_0^T \frac{\dd s}{Y_s}\bigr)^{1/2}} ,
           \frac{\int_0^T \sqrt{Y_s} \, \dd W_s}
                {\bigl(\int_0^T Y_s \, \dd s\bigr)^{1/2}} ,
           \frac{\int_0^T \sqrt{Y_s} \, \dd B_s}
                {\bigl(\int_0^T Y_s \, \dd s\bigr)^{1/2}} \Biggr) \\
   \distr
   (Z_1, -\bZ_2) \qquad
   \text{as \ $T \to \infty$,}
 \]
 see Barczy and Pap \cite[Formula (7.6)]{BarPap}.
Hence we obtain \eqref{aabb_super}.

It is known that under the condition
 \ $a \in \left[ \frac{\sigma_1^2}{2}, \infty \right)$, \ we have
 \ $\PP(\tcY_{-1/b} \in \RR_{++}) = 1$ \ (see, e.g., page 442 in Revuz and Yor
 \cite{RevYor}) and \ $\PP(\int_0^{-1/b} \tcY_s \, \dd s \in \RR_{++}) = 1$
 \ (which has been shown in the proof of Theorem 3.1 in Barczy et
 al.~\cite{BarDorLiPap}), hence \eqref{LAQJ} holds.

If \ $a \in \bigl( \frac{\sigma_1^2}{2}, \infty \bigr)$ \ and
 \ $\alpha \in \RR$ \ are fixed, then LAMN property of the subfamily will
 follow from 
 \begin{equation}\label{LAQDJsupcrit2}
  \EE\left( \exp\left\{ \bh_2^\top \bDelta^{(2)}_\btheta
                        - \frac{1}{2}
                          \bh_2^\top \bJ^{(2)}_\btheta \bh_2 \right\} \right)
  =: E_2 = 1
 \end{equation}
 for all \ $\bh_2 \in \RR^2$.
\ We have
 \begin{align*}
  E_2&=\EE\Biggl(\exp\Biggl\{\biggl(- \frac{\tcY_{-1/b}}{b}\biggr)^{-1/2}
                           \bh_2^\top
                           \begin{bmatrix}
                            \sigma_1 & \sigma_2 \varrho \\
                            0 & \sigma_2 \sqrt{1 - \varrho^2}
                           \end{bmatrix}^{-1} \!\!
                            \bZ_2
                            - \frac{1}{2}
                            \biggl(- \frac{\tcY_{-1/b}}{b}\biggr)^{-1}
                            \bh_2^\top \bS^{-1} \bh_2\Biggr\}\Biggr) \\
   &=\EE\Biggl(\EE\Biggl(\exp\Biggl\{\biggl(- \frac{\tcY_{-1/b}}{b}\biggr)^{-1/2}
                           \bh_2^\top
                           \begin{bmatrix}
                            \sigma_1 & \sigma_2 \varrho \\
                            0 & \sigma_2 \sqrt{1 - \varrho^2}
                           \end{bmatrix}^{-1} \!\!
                            \bZ_2 \\
   &\phantom{=\EE\Biggl(\EE\Biggl(\exp\Biggl\{}
                            - \frac{1}{2}
                            \biggl(- \frac{\tcY_{-1/b}}{b}\biggr)^{-1}
                            \bh_2^\top \bS^{-1} \bh_2\Biggr\}
                            \, \Bigg| \, \tcY_{-1/b} \Biggr)\Biggr)
    = 1
 \end{align*}
 by \eqref{momgennorm}, thus we conclude \eqref{LAQDJsupcrit2}.

Finally, we show that \eqref{LAQDJ} is not valid for the whole family
 \ $(\cE_T)_{T\in\RR_{++}}$ \ of statistical experiments, given in \eqref{cET},
 i.e., there exists \ $\bh \in \RR^4$, \ such that
 \begin{equation}\label{LAQDJsupcrit}
  \EE\left( \exp\left\{ \bh^\top \bDelta^{(2)}_\btheta
                        - \frac{1}{2}
                          \bh^\top \bJ^{(2)}_\btheta \bh \right\} \right)
  \ne 1 .
 \end{equation}
Indeed, using again \eqref{momgennorm}, for
 \ $\bh = (0, 1, \bzero)^\top \in \RR \times \RR \times \RR^2$, \ we have
 \begin{align*}
  &\EE\Biggl( \exp\Biggl\{ \bh^\top \bDelta_\btheta
                           - \frac{1}{2}
                             \bh^\top \bJ_\btheta \bh \Biggr\} \Biggr) \\
  &=\EE\Biggl( \exp\Biggl\{ \begin{bmatrix} 0 \\ 1 \end{bmatrix}^\top
                            \begin{bmatrix}
                             \sigma_1 & \sigma_2 \varrho \\
                             0 & \sigma_2 \sqrt{1 - \varrho^2}
                            \end{bmatrix}^{-1}
                            \begin{bmatrix}
                             \sigma_1^{-1} \tcV \\
                             Z_1
                            \end{bmatrix} \\
  &\phantom{=\EE\Biggl( \exp\Biggl\{}
                            - \frac{1}{2}
                              \biggl(\int_0^{-1/b} \tcY_u \, \dd u\biggr)
                              \begin{bmatrix} 0 \\ 1 \end{bmatrix}^\top
                              \bS^{-1}
                              \begin{bmatrix} 0 \\ 1 \end{bmatrix}
                   \Biggr\} \Biggr) \\
  &=\EE\Biggl(\exp\Biggl\{ \frac{1}{\sigma_2\sqrt{1-\varrho^2}} Z_1
                           - \frac{1}{2\sigma_2^2(1-\varrho^2)}
                             \int_0^{-1/b} \tcY_u \, \dd u\Biggr\} \Biggr) \\
  &=\EE\Biggl(\exp\Biggl\{\frac{1}
                               {\sigma_2\sqrt{1-\varrho^2}} Z_1\Biggr\} \Biggr)
    \EE\Biggl(\exp\Biggl\{- \frac{1}{2\sigma_2^2(1-\varrho^2)}
                            \int_0^{-1/b} \tcY_u \, \dd u\Biggr\} \Biggr) \\
  &=\exp\Biggl\{\frac{1}{2\sigma_2^2(1-\varrho^2)}\Biggr\}
    \EE\Biggl(\exp\Biggl\{- \frac{1}{2\sigma_2^2(1-\varrho^2)}
                            \int_0^{-1/b} \tcY_u \, \dd u\Biggr\} \Biggr)
   \ne 1 ,
 \end{align*}
 since, by Lemma 1 in Ben Alaya and Kebaier \cite{BenKeb1},
 \[
   \EE\Biggl(\exp\Biggl\{-2\mu^2\int_0^t \tcY_u \, \dd u\Biggr\}\Biggr) \\
   =\cosh(\sigma_1\mu t)^{-\frac{2a}{\sigma_1^2}}
    \exp\left\{\frac{2\mu y_0}{\sigma_1} \tanh(\sigma_1\mu t)\right\} 
 \]
 for \ $\mu, t \in \RR_+$.
\proofend

\appendix

\vspace*{5mm}

\noindent{\bf\Large Appendix}

\section{A limit theorem for continuous local martingales}

In what follows we recall a so called stable central limit theorem for
 multidimensional continuous local martingales.

\begin{Thm}{\bf (van Zanten \cite[Theorem 4.1]{Zan})}\label{THM_Zanten}
Let \ $\bigl( \Omega, \cF, (\cF_t)_{t\in\RR_+}, \PP \bigr)$ \ be a filtered
 probability space satisfying the usual conditions.
Let \ $(\bM_t)_{t\in\RR_+}$ \ be a $d$-dimensional continuous local martingale
 with respect to the filtration \ $(\cF_t)_{t\in\RR_+}$ \ such that
 \ $\PP(\bM_0 = \bzero) = 1$.
\ Suppose that there exists a function \ $\bQ : \RR_+ \to \RR^{d \times d}$
 \ such that \ $\bQ(t)$ \ is an invertible (non-random) matrix for all
 \ $t \in \RR_+$, \ $\lim_{t\to\infty} \|\bQ(t)\| = 0$ \ and
 \[
   \bQ(t) \langle \bM \rangle_t \, \bQ(t)^\top \stoch \bfeta \bfeta^\top
   \qquad \text{as \ $t \to \infty$,}
 \]
 where \ $\bfeta$ \ is a \ $d \times d$ random matrix.
Then, for each $\RR^k$-valued random vector \ $\bv$ \ defined on
 \ $(\Omega, \cF, \PP)$, \ we have
 \[
   (\bQ(t) \bM_t, \bv) \distr (\bfeta \bZ, \bv) \qquad
   \text{as \ $t \to \infty$,}
 \]
 where \ $\bZ$ \ is a \ $d$-dimensional standard normally distributed random
 vector independent of \ $(\bfeta, \bv)$.
\end{Thm}

We note that Theorem \ref{THM_Zanten} remains true if the function \ $\bQ$
 \ is defined only on an interval \ $[t_0, \infty)$ \ with some
 \ $t_0 \in \RR_{++}$.


\begin{thebibliography}{99}

\bibitem{BarDorLiPap}
\textsc{Barczy, M.}, \textsc{D\"oring, L.}, \textsc{Li, Z.} and
\textsc{Pap, G.} (2013).
On parameter estimation for critical affine processes.
\textit{Electronic Journal of Statistics}
\textbf{7} 647--696.

\bibitem{BarDorLiPap2}
\textsc{Barczy, M.}, \textsc{D\"oring, L.}, \textsc{Li, Z.} and
 \textsc{Pap, G.} (2013).
 Stationarity and ergodicity for an affine two factor model.
\textit{Advances in Applied Probability}.
\textbf{46(3)} 878--898.

\bibitem{BarPap}
\textsc{Barczy, M.} and \textsc{Pap, G.} (2013).
Maximum likelihood estimation for Heston models.
Available on ArXiv:  \texttt{http://arxiv.org/abs/1310.4783}

\bibitem{BenKeb1}
\textsc{Ben Alaya, M.} and \textsc{Kebaier, M.} (2012).
Parameter estimation for the square root diffusions: ergodic and nonergodic
 cases.
\textit{Stochastic Models}
\textbf{28(4)} 609--634.

\bibitem{BenKeb2}
\textsc{Ben Alaya, M.} and \textsc{Kebaier, M.} (2013).
Asymptotic behavior of the maximum likelihood estimator for ergodic and
 nonergodic square-root diffusions.
\textit{Stochastic Analysis and Applications}
\textbf{31(4)} 552--573.

\bibitem{CoxIngRos}
\textsc{Cox, J. C.}, \textsc{Ingersoll, J. E.} and \textsc{Ross, S. A.} (1985).
A theory of the term structure of interest rates.
\textit{Econometrica}
\textbf{53(2)} 385--407.

\bibitem{DufFilSch}
\textsc{Duffie, D.}, \textsc{Filipovi\'{c}, D.} and \textsc{Schachermayer, W.}
 (2003).
Affine processes and applications in finance.
\textit{Annals of Applied Probability}
\textbf{13} 984--1053.

\bibitem{Hes}
\textsc{Heston, S.} (1993).
A closed-form solution for options with stochastic volatilities with
 applications to bond and currency options.
\textit{The Review of Financial Studies}
\textbf{6} 327--343.

\bibitem{JSh}
\textsc{Jacod, J.} and \textsc{Shiryaev, A. N.} (2003).
\textit{Limit Theorems for Stochastic Processes}, 2nd ed.
Springer-Verlag, Berlin.

\bibitem{Jeg}
\textsc{Jeganathan, P.} (1995).
Some aspects of asymptotic theory with applications to time series models.
\textit{Econometric Theory} 
\textbf{11(5)} 818--887.

\bibitem{KarShr}
\textsc{Karatzas, I.} and \textsc{Shreve, S. E.} (1991).
\textit{Brownian Motion and Stochastic Calculus}, 2nd ed.
 Springer-Verlag.

\bibitem{KatMog}
\textsc{K\'atai, I.} and \textsc{Mogyor\'odi, J.} (1967).
Some remarks concerning the stable sequences of random variables.
\textit{Publicationes Mathematicae Debrecen}
\textbf{14} 227--238.

\bibitem{LeCamYang}
\textsc{Le Cam, L.} and \textsc{Yang, G. L.} (2000).
\textit{Asymptotics in statistics: some basic concepts},
Springer.

\bibitem{LiMa}
\textsc{Li, Z.} and \textsc{Ma, C.} (2013).
 Asymptotic properties of estimators in a stable Cox-Ingersoll-Ross model.
 Available on the ArXiv: \texttt{http://arxiv.org/abs/1301.3243}

\bibitem{LipShiI}
\textsc{Liptser, R. S.} and \textsc{Shiryaev, A. N.} (2001).
\textit{Statistics of Random Processes I. Applications}, 2nd edition.
Springer-Verlag, Berlin, Heidelberg.

\bibitem{LipShiII}
\textsc{Liptser, R. S.} and \textsc{Shiryaev, A. N.} (2001).
\textit{Statistics of Random Processes II. Applications}, 2nd edition.
Springer-Verlag, Berlin, Heidelberg.

\bibitem{Ove}
\textsc{Overbeck, L.} (1998).
Estimation for continuous branching processes.
\textit{Scandinavian Journal of Statistics}
\textbf{25(1)} 111--126.

\bibitem{RevYor}
\textsc{Revuz, D.} and \textsc{Yor, M.} (1999).
\textit{Continuous Martingales and Brownian Motion}, 3rd ed.
 Springer-Verlag Berlin Heidelberg.

\bibitem{Str}
\textsc{Stroock, D. W.} (1993).
\textit{Probability Theory, an Analytic View}.
Cambridge University Press, Cambridge.

\bibitem{Vaart}
\textsc{van der Vaart, A. W.} (1998).
\textit{Asymptotic Statistics},
Cambridge University Press.

\bibitem{Zan}
\textsc{van Zanten, H.} (2000).
A multivariate central limit theorem for continuous local martingales.
\textit{Statistics \& Probability Letters}
\textbf{50(3)} 229--235.

\end{thebibliography}
\end{document}